\documentclass[reqno]{amsart}
\usepackage{amscd}
\usepackage[dvips]{graphicx}

\newtheorem{thm}{\indent\sc Theorem}

\newenvironment{pf}{\begin{proof}[\indent \it Proof]}{\end{proof}}

\newtheorem{dfn}{\indent\sc Definition}[subsection]

\newtheorem{prp}[dfn]{\indent\sc Proposition}

\newtheorem{cor}[dfn]{\indent\sc Corollary}

\newtheorem{lma}[dfn]{\indent\sc Lemma}

\newtheorem{rem}[dfn]{\indent\it Remark}
\newenvironment{rmk}{\begin{rem}\rm}{\end{rem}}

\newcommand{\tR}{{\mathbb{R}}}

\newcommand{\tQ}{{\mathbb{Q}}}
\newcommand{\tH}{{\mathbb{H}}}
\newcommand{\tZ}{{\mathbb{Z}}}
\newcommand{\tF}{{\mathcal{F}}}
\newcommand{\Si}{\Sigma}

\newcommand{\wt}{\widetilde}
\newcommand{\wM}{{\widetilde M}}
\newcommand{\wF}{{\widetilde F}}
\newcommand{\tA}{{\mathbf A}}

\newcommand{\tM}{{\mathbf M}}
\newcommand{\tN}{{\mathbf N}}

\newcommand{\tJ}{{\mathbf J}}
\newcommand{\tI}{{\mathbf I}}

\newcommand{\lk}{\operatorname{lk}}

\newcommand{\tP}{{\mathcal P}}
\newcommand{\tT}{{\mathcal T}}
\newcommand{\Imm}{{\mathbf {Imm}}}
\newcommand{\Emb}{{\mathbf  {Emb}}}

\newcommand{\St}{\operatorname{St}}

\renewcommand{\mod}{\operatorname{mod}}

\begin{document}

\title
[3-knots with and without self intersections\dots]
{Differential 3-knots in 5-space with and without self intersections}
\author{Tobias Ekholm}
\address{Department of Mathematics, Uppsala University, S-751 06 Uppsala,
Sweden}
\email{tobias@math.uu.se}
\date{}
\begin{abstract}
Regular homotopy classes of immersions $S^3\to\tR^5$
constitute an infinite cyclic group. The classes containing
embeddings form a subgroup of index 24. The obstruction  for a generic
immersion to be regularly 
homotopic to an embedding is described in terms of geometric
invariants of its self intersection.  Geometric properties of self
intersections are used to construct two
invariants $J$ and $\St$ of generic  
immersions which are analogous to Arnold's invariants of
plane curves \cite{A1}. We prove that $J$ and $\St$ are independent first
order invariants   
and that any first order invariant is a linear combination of these.
 
As by-products, some invariants of immersions $S^3\to\tR^4$ are 
obtained. Using them, we find  restrictions on the topology of self
intersections.  
\end{abstract}
\maketitle
\section{Introduction}
\noindent
Immersions and embeddings form open subspaces of the
space of $C^r$-maps $S^k\to\tR^{k+n}$, $r\ge 1$. 
Smale \cite{S} showed that
the path components of the space of immersions 
(or, which is the same, the regular homotopy classes of immersions)
$S^k\to\tR^{k+n}$ are in one to one correspondence with the elements of
$\pi_k(V_{k+n,k})$, the $k^{th}$ homotopy group of the Stiefel
manifold of $k$-frames in $(k+n)$-space. This far-reaching result
translates problems in geometry to homotopy theory. 

Indicating the way
back to geometry, Smale suggested the following problems (\cite{S}, p.329):  
{\em ``Find explicit representatives of regular homotopy classes\dots    
What regular homotopy classes have an embedding for representative?''}.
Explicit representatives of
regular homotopy classes of immersions $S^3\to\tR^5$ are given in
Section ~\ref{secgen}.

An answer to the second problem gives information about how the inclusion
of the  space of embeddings into the space of immersions
is organized. In the case $S^3\to\tR^5$, the group $\pi_3(V_{5,3})$,
enumerating regular homotopy classes is infinite cyclic and the answer 
to the second problem was found by Hughes and Melvin \cite{HM}. They
proved that exactly every $24^{\rm th}$ regular homotopy class
have an embedding for representative. In Theorem \ref{thmimem} below we 
describe the obstruction for a generic immersion $S^3\to\tR^5$ to be
regularly homotopic to an embedding in
terms of geometric invariants of its self intersection.

In dimensions where there are embeddings in different regular homotopy
classes, it is impossible to express the regular homotopy class of a
generic immersion in terms of its self intersection. 
This is in contrast to many other cases where this is possible: For
example, in the cases $S^k\to\tR^{2k}$, $k\ge 2$ the regular homotopy
class is determined by the algebraic number of 
self intersection points (modulo 2 if $k$ is odd), see \cite{S}, and
in the cases  
$S^k\to\tR^{2k-r}$, $r=1,2$ self intersection formulas for regular
homotopy are given in \cite{E2} and \cite{E3}.  

Generic immersions have simple self intersections. For example,
in the case $S^1\to\tR^3$ a generic immersion has empty self
intersection and thus, generic immersions are embeddings. 
Form this point of
view, the analogue of classical knot theory in
other 
dimensions is the study of path components of the space of generic
immersions and we may think of generic immersions as knots with self
intersections. 

The space of generic immersions 
is dense in the corresponding space of immersions 
and its complement is a stratified hypersurface. Using the
stratification of this complement,
Vassiliev \cite{A1} introduced the notion of finite order invariants in 
classical knot theory. There are natural analogies of this notion
for invariants of generic immersions in other dimensions (see Section
~\ref{1stord}). Arnold
\cite{A2} found first 
order invariants of generic plane curves ($S^1\to\tR^2$).  Theorem
~\ref{thmJSt} below shows that, up 
to first order, the space of
generic immersions $S^3\to\tR^5$ is similar to the space
of generic plane curves.

The composition of an immersions $S^3\to\tR^4$ and the inclusion
$i\colon\tR^4\to\tR^5$ is an immersion into 5-space. Proposition
~\ref{4ss5} shows that two 
immersions of $S^3$ into 
4-space are regularly homotopic in 5-space, after composing them with the
inclusion, if and only if one of them is regularly homotopic in 4-space
to the connected sum of the other one and a finite number of
immersions regularly homotopic to the composition of the
standard embedding and a reflection in a hyperplane in $\tR^4$. 

Theorems
~\ref{gi41}, ~\ref{gi42}, and ~\ref{gi43} give information about the
self intersections of a generic immersions $g\colon S^3\to\tR^4$ and
its relation to the self intersection of a generic immersion $f\colon
S^3\to\tR^5$ regularly homotopic to $i\circ g$.

\section{Main results}
\noindent
In this section, the main theorems of the paper are formulated.

\subsection{Embeddings in the space of immersions}\label{emimsp}
\noindent
Regular homotopy classes of immersions $S^3\to\tR^5$ are known to form
an infinite 
cyclic group $\Imm$ under connected sum. The classes that contain
embeddings form a subgroup $\Emb\subset\Imm$. We now state Theorem
~\ref{thmimem} describing the extension $\Emb\to\Imm$ algebraically
and devote the rest of this section to identify the 
homomorphisms involved there in topological terms.
\begin{thm}\label{thmimem}
The following diagram of Abelian groups has exact rows, commutes and
the vertical arrows are isomorphisms
$$
\begin{CD}
0 @>>> \Emb @>>> \Imm @>{\lambda\oplus\beta}>>\tZ_3\oplus\tZ_8 @>>> 0\\
@VVV @V{\sigma}VV  @VV{\Omega}V @VVV @VVV\\
0@>>> \tZ @>{\times24}>> \tZ @>>> \tZ_{24} @>>> 0
\end{CD}\quad.
$$
\end{thm}
\noindent
Theorem ~\ref{thmimem} is proved in Section ~\ref{pfimem}. 

{\em The homomorphism $\sigma$}:
Recall that if
$f\colon S^3\to \tR^5$ is an embedding then there exists a compact orientable
4-dimensional manifold $V^4\subset\tR^5$ with $\partial
V^4=f(S^3)$. We call such a manifold a {\em Seifert-surface} of $f$.
Its signature $\sigma(V^4)$ is divisible by 16 and is known to depend
only on $f$.  
{\em For $\xi\in\Emb$, define $\sigma(\xi)\in\tZ$ as
$\frac{\sigma(V)}{16}$, where $V$ is a  Seifert-surface of an
embedding representing $\xi$.}  
It is proved in \cite{HM} that $\sigma$ induces an
isomorphism $\Emb\to\tZ$ and that $\Emb\subset\Imm$ is a subgroup of
index $24$.

{\em The homomorphism $\Omega$}:
Any immersion $f\colon S^3\to\tR^5$ is determined up to regular
homotopy by its Smale
invariant $\Omega(f)\in\pi_3(V_{5,3})\cong\tZ$ (see \cite{S} or
Definition ~\ref{dfnSmin}). 
{\em For $\xi\in\Imm$,  define $\Omega(\xi)\in\tZ$ as $\Omega(g)$,
where $g$ is an immersion representing $\xi$.} 

{\em The homomorphism $\lambda$}:
If $f\colon S^3\to\tR^5$ is a generic immersion then its self
intersection $M_f\subset\tR^5$ is a closed 1-dimensional
manifold. Orientations of $S^3$ and
$\tR^5$ induce an orientation of $M_f$. 
We can push $M_f$ off the
image of $f$ (see Section ~\ref{seclk}). 
Let  $\lk(f)$ denote the linking number of the perturbed $M_f$ and
$f(S^3)$ in $\tR^5$ .
{\em For $\xi\in\Imm$, define $\lambda(\xi)\in\tZ_3$ as $\lk(g)$
modulo 3, where $g$ is a generic immersion representing $\xi$.}

{\em The homomorphism $\beta$}:
The normal bundle of an immersion $f_0\colon S^3\to\tR^5$ is
2-dimensional orientable and therefore trivial. Hence, $f_0$ admits a
normal vector field. This 
implies that $f_0$ is regularly homotopic to a (generic)
immersion $f_1\colon S^3\to\tR^4$, composed with the inclusion
$\tR^4\to\tR^5$ (see \cite{Hi}). 
Resolving the self intersection of a generic immersion  
$f\colon S^3\to\tR^4$, we obtain a smooth surface $F_f$ (Lemma
~\ref{FgCg}) and the immersion $f$ induces a pin (i.e. $Pin^-$) 
structure on $F_f$ (Section ~\ref{secbeta}).  There is a
one to one correspondence between pin structures on a surface $F$ and
$\tZ_4$-quadratic functions $q$ on its first homology $H_1(F;\tZ_2)$. Pin
structures on a surface are classified up to cobordism by the Brown
invariant $\beta(q)\in\tZ_8$ of the corresponding function.
Let $\beta(f)\in\tZ_8$ denote the Brown invariant of the quadratic function
corresponding to the pin structure induced by $f$ on $F_f$. 
{\em For $\xi\in\Imm$, define
$\beta(\xi)\in\tZ_8$ as $\beta(g)$, where $g\colon
S^3\to\tR^4\subset\tR^5$ is a generic immersion representing $\xi$.}
The author does not know how to calculate $\beta(f)$ for a generic
immersion $f\colon S^3\to\tR^5$ in terms of the geometry of its self
intersection, without pushing it down to $\tR^4$. However, 
$\beta(f)$ modulo 4, can be calculated in that way (see Section
~\ref{secrelinv}). 

\subsection{First order invariants of generic immersions}
\noindent
In generic one-parameter families of immersions there are isolated
instances of non-generic immersions.
In generic one-parameter families of immersions
$S^3\to\tR^5$ such instances are immersions with one self
tangency or one triple point. The same is true for generic one-parameter
families of plane curves where one can distinguish two local types of
self tangencies:  
direct (the tangent vectors point in the same direction) and reverse
(the tangent vectors point in opposite directions). 

An {\em invariant of generic immersions} is a function which is
constant on the path 
components of the space of generic immersions. Such a function may
change when we pass through instances of non-generic immersions. 

Arnold \cite{A2} found three independent first order invariants of
generic plane curves:  
The invariant  $J^+$ which changes
at direct self tangency instances and does not change under other
moves,  the invariant $J^-$ which changes under reverse self tangency
moves and does not change under other moves, and the invariant
Strangeness $\St$ 
which changes under triple point moves and does not change under other
moves. 

For a generic immersion $f\colon S^3\to\tR^5$, define $J(f)$ to be the
number of self intersection components of $f$ and define
$L(f)=\frac13(\lk(f)-{\wt \lambda}(f))$, ${\wt
\lambda}\in\{0,1,2\}\subset\tZ$  is a lifting of $\lambda\in\tZ_3$
(see Section ~\ref{emimsp}).
\begin{thm}\label{thmJSt}
The invariant $J$ changes by $\pm 1$ under self tangency moves and does
not change under triple point moves. The invariant $L$ changes by $\pm1$
under triple point moves and does not change under self tangency
moves. The invariants $J$ and $L$ are first order invariants. 
Moreover, if $v$ is any first order invariant then the restriction of
$v|U$, where $U$ is a path component of the space of immersions, is a
linear combination of $J|U$ and $L|U$.   
\end{thm}
\noindent
Theorem ~\ref{thmJSt} is proved in Section ~\ref{pfJSt}. Although
there are two local types of selftangencies in the case $S^3\to\tR^5$
(see Proposition ~\ref{verdef51}), in contrast
to the case  
of plane curves, $J$ can not be splitted into more refined first order
invariants.  

Arnold defined the invariant $\St$ in such a way that it is additive
under connected summation of plane curves and showed that this
property together with the changes of $\St$ under local moves and 
its is orientation independence
completely characterizes the invariant (up to a constant).

The invariant $L$ is neither additive under connected summation nor
symmetric with respect to orientation. To get the analogue of Arnold's
strangeness of plane curves for immersions $f\colon S^3\to\tR^5$ we define
$\St(f)=\frac13(\lk(f)+\Omega(f))$. Then $\St$ is additive under
connected sum, changes exactly as $L$ under local moves, changes sign
under composing immersions with an orientation reversing
diffeomorphism of $S^3$, and is
completely characterized by these properties up to a constant
(see Proposition ~\ref{Strange}).

The reason for stating Theorem ~\ref{thmJSt} in terms of $L$ instead
of $\St$ (which works equally well) is that $L(f)$ in contrast to
$\St(f)$ can be calculated in terms of the self intersection of $f$.

\subsection{Immersions into 4-space} 
\noindent
The self intersection of a generic immersion\linebreak
$f\colon S^3\to\tR^4$ consists of 2-dimensional sheets of
double points, 1-dimensional curves of triple points and isolated
quadruple points (see Definition ~\ref{dfngi4}). Resolving the quadruple and
triple points turns the self intersection into a smooth surface $F_f$
(Lemma ~\ref{FgCg}). 

The following three theorems are proved 
in Section ~\ref{secrelinv}.
\begin{thm}\label{gi41}
A generic immersion $g\colon S^3\to\tR^4$ has an odd number of
quadruple points if and only if its (resolved) self intersection
surface $F_g$ has odd Euler characteristic.
\end{thm}

As mentioned (Section ~\ref{emimsp}), any immersion $S^3\to\tR^5$ is
regularly homotopic to a composition of an immersion $S^3\to\tR^4$ and
the inclusion $\tR^4\to\tR^5$. 
\begin{thm}\label{gi42}
Let $f\colon
S^3\to\tR^5$ be a generic immersion and let $g\colon 
S^3\to\tR^4\subset\tR^5$ be a generic immersion regularly homotopic to
$f$. Then
$g$
has an odd number of quadruple points if and only if $f$ has an odd
number of self intersection components with connected preimage.
\end{thm}
\noindent
We define a $\tZ_4$-valued invariant $\tau$ of generic immersions
$S^3\to\tR^5$ (see Section ~\ref{tau}). If $f$ is a generic immersion
then $\tau(f)$ is divisible by 2 if and only if $f$ has an odd
number of self intersection components with connected preimage.  
\begin{thm}\label{gi43}
Let $f\colon
S^3\to\tR^5$ be a generic immersion and let $g\colon 
S^3\to\tR^4\subset\tR^5$ be a generic immersion regularly homotopic to
$f$.
If $\tau(f)\ne 0$ then the (resolved) self intersection surface $F_g$
of $g$ is nonorientable. 
\end{thm}

\section{The Smale invariant and Stiefel manifolds}
\noindent
In this section, some properties of
the Smale invariant and 
the groups where it takes values are collected for later reference. 
Via the Smale invariant the set of
regular homotopy classes of 
immersions is endowed with the structure of an Abelian group. We
describe the group operations in topological terms.  Most of the
results presented here are known. Proofs are provided where there were
hard to find references. 

\subsection{Homotopy groups of two Stiefel manifolds}\label{sechomth}
\noindent
The Stiefel manifold $V_{4,3}$ is homotopy equivalent to $SO(4)$. Indeed,
any orthonormal 3-frame in 4-space can be uniquely completed to a
positively oriented $4$-frame.

Let $SO(4)\xrightarrow{p} S^3$ be the fibration with fiber $SO(3)$.
Consider $S^3$ as the set of unit quaternions in $\tH\cong\tR^4$,
where we identify the vectors $\partial_1,\dots,\partial_4$ in the
standard base of $\tR^4$ with $1,i,j,k\in\tH$.  For $x,y\in\tR^4$
let $x\cdot y\in\tR^4$ denote  
quaternionic product of $x$ and $y$. The map $\sigma\colon S^3\to
SO(4)$, $\sigma(x)y=x\cdot y$, is a section of  
$SO(4)\xrightarrow{p}S^3$. Thus, $SO(4)$ is diffeomorphic to $S^3\times
SO(3)$.  

Let $\varrho\colon S^3\to SO(3)$ be the map $\varrho(x)u=x\cdot u\cdot
x^{-1}$, where $u$ is a pure quaternion and $\tR^3$ is identified with
the set of pure quaternions. (That is, the span of 
the vectors $i,j,k\in\tH$.) Let $\rho\colon S^3\to SO(4)$ be
$\varrho$ 
composed with the inclusion of the fiber over $(1,0,0,0)\in
S^3$. The following is immediate.
\begin{lma}\label{pi3V43}
$$
\pi_3(SO(4))\cong\tZ[\sigma]\oplus\tZ[\rho].
$$\qed
\end{lma}

Let $SO(5)\xrightarrow{p} V_{5,3}$ be the fibration with fiber $SO(2)$
that maps an orthonormal $5\times 5$-matrix to the
3-frame consisting of its
first three column vectors. Similarly, let $SO(5)\xrightarrow{r} S^4$
be the fibration with fiber $SO(4)$.
\begin{lma}\label{pi3V53}
The homomorphism $p_\ast\colon\pi_3(SO(5))\to\pi_3(V_{5,3})$ is an
isomorphism. The homomorphism
$r_\ast\colon\pi_3(SO(4))\to\pi_3(SO(5))$ is an epimorphism with
kernel $N$, where $N$ is the subgroup generated by $[\sigma]-2[\rho]$.
\end{lma}
\begin{pf}
The first statement follows by inspecting the homotopy sequence of the
fibration. For the second, see \cite{Hu}, Chapter 8, Proposition 12.11.
\end{pf}
Using Lemma ~\ref{pi3V53}, we make the identifications
$$
\pi_3(V_{5,3})=\pi_3(SO(5))=\pi_3(SO(4))/N.
$$

\subsection{The Smale invariant of an immersed 3-sphere}
\noindent
We define the Smale invariant of an immersion $f\colon S^3\to\tR^n$:
Consider $S^3$ as the unit sphere in $\tR^4$ and
let $s\colon S^3\to\tR^4\subset\tR^n$, denote the standard
embedding. Fix a 
disk $D^3\subset S^3$ containing the south pole and a framing $X$ of
$S^3-D^3$. Using regular homotopy, deform $f$ so that $f|D^3=s|D^3$. 

Choose a diffeomorphism $r\colon H_+\to S^3-D^3$ of degree $+1$, where
$H_+$ is the hemisphere $\{x_0\ge 0\}$ in the unit sphere
$\{x_0^2+x_1^2+x_2^2+x_3^2=1\}$ in $\tR^4$. Let $x\mapsto x^\ast$ be
the map $(x_0,x_1,x_2,x_3)\mapsto (-x_0,x_1,x_2,x_3)$ of the unit
sphere in $\tR^4$. Define $\phi_s^f\colon S^3\to V_{n,3}$,
$$
\phi_s^f(x)=
\begin{cases}
df(X(r(x))) & \text{for }x\in H_+,\\
ds(X(r(x^\ast))) & \text{for }x\in H_-,
\end{cases}
$$
where $H_-$ is the hemisphere $\{x_0\le0\}$.
\begin{dfn}\label{dfnSmin}
The Smale invariant $\Omega(f)$ of $f$ is
$$
\Omega(f)=[\phi^f_s]\in\pi_3(V_{n,3}).
$$
\end{dfn}
Smale, \cite{S} showed that $\Omega$ gives a bijection between the
regular homotopy classes of immersions $S^3\to\tR^n$ and the elements
of $\pi_3(V_{n,3})$. 

\subsection{Calculating Smale invariants in 4- and 5-space}\label{calc}
\noindent
In computations we will not use Definition ~\ref{dfnSmin}
literally. We use a slightly different approach: 
Consider 
$$
S^3=\{x\in\tR^4=\tH : x_0^2+x_1^2+x_2^2+x_3^2=1\},
$$ 
where $\tH$ denotes the quaternions.
The tangent
space of $S^3$ at $(1,0,0,0)$ is the span of the vectors $i,j,k$. We
trivialize $TS^3$ using the {\em quaternion framing:}
$$
Q(x)=(x\cdot i,x\cdot j,x\cdot k)\in T_xS^3,\text{ for $x\in S^3$.}
$$

Let $f\colon S^3\to \tR^n$ be an immersion. Then there is an induced
map $\Phi_n^f\colon S^3\to V_{n,3}$,  $\Phi_n^f(x)=df(Q(x))$. 

Assume that $n=4$. Then $V_{n,3}=V_{4,3}=SO(4)$. Thus, we get a map
$\Phi_4^f\colon S^3\to SO(4)$.  
\begin{lma}\label{calcSmin4}
Let $f\colon S^3\to\tR^4$ be an immersion then
$$
\Omega(f)=[\Phi_4^f]-[\Phi_4^s]\in\pi_3(SO(4)),
$$  
and $[\Phi_4^s]=[\sigma]\in\pi_3(SO(4))$.\qed
\end{lma}

Assume that $n=5$. Consider the
fibration $SO(5)\to V_{5,3}$ as in 
Lemma ~\ref{pi3V53}.
Since the normal bundle of $f$ 
is orientable 2-dimensional, it is trivial and we can lift $\Phi_5^f$
to $\Theta^f\colon S^3\to SO(5)$ and $[\Theta^f]\in\pi_3(SO(5))$ is
independent of this lifting.    
\begin{lma}\label{calcSmin5}
Let $f\colon S^3\to\tR^5$ be an immersion then
$$
\Omega(f)=[\Theta^f]-[\Theta^s]\in\pi_3(SO(5))=\pi_3(V_{5,3}),
$$ 
and $[\Theta^s]=[\sigma]+N\in\pi_3(SO(5))$.\qed
\end{lma}

Let $i\colon\tR^4\to\tR^5$ denote the inclusion.
\begin{lma}\label{Smin4in5}
If $f\colon S^3\to\tR^4$ is an immersion then
$$
\Omega(i\circ f)=\Omega(f)+N\in\pi_3(SO(5))=\pi_3(V_{5,3}).
$$
\end{lma}
\begin{pf}
Complete the framing of $f$
with a vector in the fifth direction and combine Lemma
~\ref{calcSmin4} and Lemma ~\ref{calcSmin5}.
\end{pf}

\subsection{The immersion group}
\noindent
Let $\Imm$ denote the infinite cyclic group of regular homotopy
classes of immersions $S^3\to\tR^5$.
The Smale invariant gives an isomorphism 
$\Omega\colon\Imm\to\pi_3(V_{5,3})$.

First, we consider addition in $\Imm$ and other groups of regular
homotopy classes of immersions $S^3\to\tR^n$:

Given two immersions $f,g\colon S^3\to\tR^n$ we define an immersion
$f\star g$ as follows:
Consider $S^3\subset\tR^4$ with coordinates $x=(x_0,x_1,x_2,x_3)$ as the
subset characterized by $\sum x_i^2=1$. Let $a=(1,0,0,0)$ and
$a'=(-1,0,0,0)$, respectively. Choose frames $u_1,\dots,u_n$ at $f(a)$
and $v_1,\dots,v_n$ at $g(a')$ that agree with the orientation of
$\tR^n$ and such that $u_1,u_2,u_3$ ($v_1,v_2,v_3$) are
tangent to $f(S^3)$ (to $g(S^3)$) at $f(a)$ (at $g(a')$). We can
assume, possibly after moving $g(S^3)$, that $u_i=v_i$, $i=1,2,3$ and
that $g(a')=f(a)+u_n$. Moreover, we can deform the maps so that
\begin{align*}
f(x) &= f(a)+\sum_{i=1}^3x_iu_i\quad\text{for}\quad 1-\epsilon\le
x_0\le 1,\\
g(x) &= g(a')+\sum_{i=1}^3x_iv_i\quad\text{for}\quad -1\le
x_0\le -1+\epsilon.
\end{align*}  
The immersion $f\star g$ is now obtained by running a tube from $f(a)$
to $g(a')$ with axis $f(a)+tu_n$. Details can be found in Kervaire
\cite{K}, Section 2, where the following is proved.
\begin{lma}\label{sum}
$$
\Omega(f\star g)=\Omega(f)+\Omega(g).
$$\qed
\end{lma}
We call the immersion $f\star g$ {\em the connected sum of $f$ and $g$}.

Secondly, we consider inversion in $\Imm$: 

Given an immersion $f\colon S^3\to\tR^5$ we define the immersion
$\widehat f$ as follows: Let $\widehat f=f\circ r$, where
$r\colon S^3\to S^3$ is the restriction of a reflection through a
hyperplane in $\tR^4$.
\begin{lma}\label{invrs}
$$
\Omega({\widehat f} )=-\Omega(f)\in\pi_3(V_{5,3}).
$$
\end{lma}
\begin{pf}
Let the disk $U\subset S^3\subset\tR^4$ in the definition of the Smale
invariant be the hemisphere $\{x_0\ge 0\}$ and let $r\colon S^3\to
S^3$ be the map 
$$
r(x_0,x_1,x_2,x_3)=(x_0,-x_1,x_2,x_3).
$$ 
Let $R\colon\tR^5\to\tR^5$ be the map
$$
R(y_1,y_2,y_3,y_4,y_5)=(y_1,-y_2,y_3,y_4,-y_5).
$$ 
Then $R\circ s\circ r=s$ and we can use the map 
$\phi^{R\circ f\circ r}_{R\circ s\circ r}\colon S^3\to V_{3,5}$ to
compute $\Omega({\widehat f} )$. It is straightforward to check that
this map is homotopic to $\phi_s^f\circ r$. Hence, 
$$
\Omega({\widehat f} )=(\phi^{R\circ f\circ r}_{R\circ s\circ r})_\ast[S^3]=
(\phi_s^f\circ r)_\ast[S^3]=
(\phi_s^f)_\ast\left(-[S^3]\right)=-\Omega(f).
$$  
\end{pf}

\section{Embeddings considered as immersions}
\noindent
In this section we present the classification of embeddings
$S^3\to\tR^5$ up to regular homotopy.

\subsection{Embeddings up to regular homotopy and signature} 
\noindent
Let $f,g\colon S^3\to\tR^5$ be embeddings. Then $f\star g$ is
regularly homotopic to an embedding and ${\widehat f}$ is an
embedding. Thus, the regular homotopy classes that contain embeddings
from a subgroup of $\Imm$. We denote this subgroup $\Emb$.

We get a classification of embeddings up to regular homotopy as follows:

Given an embedding $f\colon S^3\to\tR^5$, we can find a compact
connected orientable manifold $V^4$ embedded into $\tR^5$ and such
that $\partial V^4=f(S^3)$. We call such a manifold a 
{\em Seifert-surface of $f$}. 
The orientation of $f(S^3)$ induces an orientation of $V^4$.  
Filling the $S^3$ on the boundary of $V^4$ with a 4-disk we
get a closed connected oriented 4-manifold $W^4$. The cohomology
sequence of the pair $(W^4,V^4)$ shows that the inclusion induces an
isomorphism 
$H^2(W^4;\tZ_2)\to H^2(V^4;\tZ_2)$. The normal bundle of $V^4$ is a trivial
1-dimensional bundle. Hence, $w_2(TV)=0$ and therefore $w_2(TW)=0$,
where $w_2$ is the second Stiefel-Whitney class. Thus, $W$ is a spin
manifold. By Rokhlin's theorem (see Milnor and Kervaire, \cite{MK}),
{\em the signature $\sigma(W^4)$ of $W^4$} is divisible by 16.  

Define $\sigma(f)=\frac{\sigma(W)}{16}\in\tZ$. (This definition agrees
with that given in Section ~\ref{emimsp}.)

\begin{prp}\label{sgntr}
For $\xi\in\Emb$, let $\sigma(\xi)=\sigma(f)$, where $f$ is an
embedding representing $\xi$. Then
$$
\sigma\colon\Emb\to\tZ 
$$
is an isomorphism.
\end{prp}
\begin{pf}
See \cite{HM}.
\end{pf}
\begin{prp}\label{24}
The subgroup $\Emb\subset\Imm$ has index $24$.
\end{prp}
\begin{pf}
See \cite{HM}.
\end{pf}

\section{Spaces of immersions}
\noindent
In this section, generic immersions $S^3\to\tR^n$, $n=4,5$ and
their self intersections are studied. The space of immersions
$S^3\to\tR^n$ as described 
in the Introduction, will be  
denoted $\tF_n$. It is an infinite dimensional manifold. The set of
non-generic immersions in $\tF_n$ is a stratified hypersurface
$\Si_n$. We describe its strata of codimension one, for $n=4,5$, and
of codimension two, for $n=5$.   

\subsection{Generic immersions and their self intersections}\label{genim}
\noindent
\begin{dfn}\label{dfngi4}
An immersion $f\in\tF_4$ is {\em generic} if it satisfies
the following conditions:
\begin{itemize}
\item[{\bf g1}] For any $w\in\tR^5$, $f^{-1}(w)$ contains at most four
points. 
\item[{\bf g2}] If $f(x_1)=\dots=f(x_j)=w$, $2\le j\le4$ for $x_i\ne
x_j\in S^3$ if $i\ne j$, then 
$df(T_{x_i}S^3)+\bigcap_{j\ne i}df(T_{x_j}S^3)=T_w\tR^4$.  
\end{itemize}
\end{dfn}

\begin{dfn}
An immersion $f\in\tF_5$ is {\em generic} if it satisfies
the following conditions:
\begin{itemize}
\item[{\bf G1}] For any $w\in\tR^5$, $f^{-1}(w)$ contains at most two
points. 
\item[{\bf G2}] If $f(x)=f(y)=w$, for $x\ne y\in S^3$, then
$df(T_xS^3)+df(T_yS^3)=T_w\tR^5$.  
\end{itemize}
\end{dfn}

If $f\colon X\to Y$ is an immersion of manifolds then the {\em self
intersection} 
of $f$ is the subset of points $y\in Y$ such that $f^{-1}(y)$
contains more than one point. We denote it $M_f$. We denote its
preimage ${\wt M}_f$. That is, ${\wt M}_f=f^{-1}(M_f)\subset X$. 

\begin{lma}
Let $f\colon S^3\to\tR^5$ be a generic immersion. Then $M_f$ and ${\wt
M}_f$ are closed 1-manifolds and $f\colon \wM_f\to M_f$ is a double
cover. Moreover, there is an induced orientation on $M_f$. 
\end{lma}   
\begin{pf}
The first statement follows from {\bf G2}.
If we order the oriented sheets coming together along $M_f$ then
there is a standard way to assign an orientation to $M_f$. Since the
codimension of ${\wt M}_f$ is even this orientation is independent of
the ordering.
\end{pf}
\begin{lma}
Let $g\colon S^3\to\tR^4$ be a generic immersion. Then $M_g$ and ${\wt
M}_g$ are
2-dimensional stratified spaces, 
$$
M_g=M_g^0\cup M_g^1\cup M_g^2\quad\text{and}\quad
{\wt M}_g={\wt M}_g^0\cup {\wt M}_g^1\cup {\wt M}_g^2,
$$
where $M_g^j$ and ${\wt M}_g^j$ are smooth manifolds of dimension $j$,
$j=0,1,2$. The strata $M_g^{j}$ is the set of $j$-tuple points, ${\wt
M}_g^j=g^{-1}(M_g^j)$, and $g|{\wt M}_g^j$ is a $j$-fold covering. 
\end{lma}
\begin{pf}
Immediate from {\bf g1} and {\bf g2}.
\end{pf}

The next lemma shows how to resolve the self intersection of a generic
immersion $S^3\to\tR^4$.
\begin{lma}\label{FgCg}
Let $g\colon S^3\to\tR^4$ be a generic immersion. There exist
closed surfaces ${\wt F}_g$ and $F_g$ (and closed 1-manifolds ${\wt C}_g$
and $C_g$), unique up to diffeomorphisms, and
immersions $s\colon{\wt F}_g\to S^3$ and $t\colon F_g\to\tR^4$ 
($\sigma\colon{\wt C_g}\to S^3$ and $\tau\colon C_g\to\tR^4$)
such that 
$$
\begin{CD}
\wF_g @>s>> \wM_g\subset S^3\\
@VpVV  @VVgV\\
F_g @>t>> M_g\subset\tR^5
\end{CD}\quad,\quad\quad
\begin{CD}
{\wt C}_g @>{\sigma}>> \wM_g^0\cup \wM_g^1\subset S^3\\
@V{\pi}VV  @VVgV\\
C_g @>{\tau}>> M_g^0\cup\wM_g^1\subset\tR^5
\end{CD}\quad,
$$
commutes. The maps $s$ and $t$ are surjective, have multiple points
only along $\wM_g^0\cup\wM_g^1$ and $M_g^0\cup M_g^1$, respectively and 
$p$ is the orientation double cover. The maps  
$\sigma$ and $\tau$ are surjective, have multiple points
only along $\wM_g^0$ and $M_g^0$, respectively and 
$\pi$ is a $3$-fold cover.  
\end{lma}
\begin{pf}
This is immediate from the local pictures: Close to a $j$-tuple point
$M_g$  is the intersection of $j$ 3-planes in general position in
4-space.  
\end{pf}
We call $F_g$ the {\em resolved self intersection surface} of
$g$. We use the notation $F^j_g=t^{-1}(M_g^j)$ and
$\wF_g^j=s^{-1}(\wM_g^j)$, $j=0,1$.  

\subsection{The discriminant hypersurface and its stratum of
codimension one} 
\noindent
The jet transversality theorem implies that the set of generic
immersions is an open dense subset of $\tF_n$, $n=4,5$. 
Its complement $\Si_n\subset\tF_n$ will be called {\em the discriminant
hypersurface}. The discriminant hypersurface is stratified,
$\Si_n=\Si^1_n\cup \Si^2_n\cup\dots\cup \Si^\infty_n$, where each
stratum $\Si^k_n$, 
$k<\infty$ is a smooth submanifold of codimension $k$ of $\tF_n$ and
$\Si^j_n$ is contained in the closure of $\Si^1_n$ for every $j$.

The following two propositions follow by applying the jet
transversality theorem  to 1-parameter families of immersions.
\begin{prp}\label{codim14}
The codimension one stratum $\Si^1_4\subset\Si_4$ is the set of all
immersions $f\colon S^3\to\tR^4$ such that
\begin{itemize}
\item[{\rm (a)}]{\bf g1} and  {\bf g2} holds
except at one double point 
$w=f(x)=f(y)\in\tR^4$, $x\ne y\in
S^3$  where 
$$
\dim \left(df(T_xS^3) + df (T_yS^3)\right)=3
$$
or,
\item[{\rm (b)}]{\bf g1} and  {\bf g2} holds
except at one $j$-tuple point, $3\le j\le 5$ 
$w=f(x_1)=\dots=f(x_j)\in\tR^4$, $x_i\ne x_k\in
S^3$, if $i\ne k$  where 
$$
\dim \left(df(T_{x_i}S^3)+\bigcap_{r\ne i\,r\ne k}df(T_{x_r}S^3)\right)=4,
$$ 
for $k\ne i$ but 
$$
\dim \left(df(T_{x_i}S^3)+\bigcap_{r\ne i}df(T_{x_r}S^3)\right)=3.
$$ 
\end{itemize}\qed
\end{prp}
\begin{prp}\label{codim15}
The codimension one stratum $\Si^1_5\subset\Si_5$ is the set of all
immersions $f\colon S^3\to\tR^5$ such that either
\begin{itemize}
\item[{(\rm a)}] {\bf G1} holds, $f$ has double points and {\bf G2}
holds, except at one double point $w=f(x)=f(y)\in\tR^5$, $x\ne y\in
S^3$, where $\dim \left(df(T_xS^3)+df(T_yS^3)\right)=4$ or,
\item[{(\rm b)}] {\bf G2} holds and {\bf G1} holds except at one triple point
$w\in\tR^5$ with $w=f(x_1)=f(x_2)=f(x_3)$, where $x_i$, $1\le i\le 3$ are three
distinct points in $S^3$ and 
$$
\dim \left(df(T_{x_i}S^3)+\bigcap_{j\ne i}df(T_{x_j}S^3)\right)=4.
$$
\end{itemize}\qed
\end{prp}

If {(\rm a)} above holds for an immersion $f$,
we say that the exceptional double point $w$ is a {\em self tangency point}
of $f$. 

\subsection{Coordinate expressions and versal deformations}\label{vdf1}
\noindent
Recall that
a {\em deformation of a map $f_0$} is a 1-parameter family of maps
$f_\lambda$, parameterized by $\lambda\in U$, where $U$ is a
neighborhood of $0\in\tR^m$. 
A deformation $f_\lambda$ of a map $f_0$ is called {\em versal} if every
deformation of $f_0$ is equivalent (up to left-right actions of
diffeomorphisms) to one induced from $f_\lambda$. 

Below, $f_0\colon S^3\to\tR^n$, $n=4,5$ will be an immersion in
$\Si_n^1$ with its exceptional multiple point at $0\in\tR^n$,
$x,y,z,w,u$ will denote coordinates in small 3-balls centered at the
preimages of $0$ and $f_t$, $t\in U\subset\tR$ will be a versal
deformation of $f_0$.  
Such a deformation can be assumed to be constant in $t$ outside of the
coordinate balls. The proofs of the statements in this section and the
next one are discussed in Section ~\ref{discpf}.
\begin{prp}\label{verdef41}
Let $f_0\colon S^3\to\tR^4$ be an immersion in $\Si_4^1$ with
exceptional multiple point $0$ and let $f_t$ be a versal
deformation. Locally, at $0$, up to 
choice of coordinates around the preimages of 
$0$ and in $\tR^4$, 
\begin{itemize}
\item[{\rm (a)}] if $0$ is double point then $f_t$ is of the form
\begin{align*}
f_t(x)&=(x_1,x_2,x_3,0),\\
f_t(y)&=(y_1,y_2,y_3,y_1^2+y_2^2+\epsilon y_3^2+t),
\end{align*}
where $\epsilon=\pm1$, or
\item[{\rm (b)}] if $0$ is a triple point then $f_t$ is of the form
\begin{align*}
f_t(x)&=(x_1,x_2,x_3,0),\\
f_t(y)&=(y_1,y_2,0,y_3),\\
f_t(z)&=(z_1,z_2,z_3,z_1^2+\epsilon z_2^2+t-z_3),
\end{align*}
where $\epsilon=\pm1$, or
\item[{\rm (c)}] if $0$ is a quadruple point then $f_t$ is of the form
\begin{align*}
f_t(x)&=(x_1,x_2,x_3,0),\\
f_t(y)&=(y_1,y_2,0,y_3),\\
f_t(z)&=(z_1,0,z_2,z_3),\\
f_t(w)&=(w_1,w_2,w_3,w_1^2+t-w_2-w_3),
\end{align*}
or
\item[{\rm (d)}] if $0$ is a quintuple point then $f_t$ is of the form
\begin{align*}
f_t(x)&=(x_1,x_2,x_3,0),\\
f_t(y)&=(y_1,y_2,0,y_3),\\
f_t(z)&=(z_1,0,z_2,z_3),\\
f_t(w)&=(0,w_1,w_2,w_3),\\
f_t(u)&=(u_1,u_2,u_3,t-u_1-u_2-u_3).
\end{align*}
\end{itemize}\qed
\end{prp}

\begin{prp}\label{verdef51}
Let $f_0\colon S^3\to\tR^5$ be an immersion in $\Si_5^1$ with
exceptional multiple point $0$ and let $f_t$ be a versal
deformation. Locally, at $0$, up to 
choice of coordinates around the preimages of 
$0$ and in $\tR^5$, 
\begin{itemize}
\item[{\rm (a)}] if $0$ is double point then $f_t$ is of the form
\begin{align*}
f_t(x)&=(x_1,x_2,0,x_3,0).\\
f_t(y)&=(y_1,y_2,y_1^2+\epsilon y_2^2+t,0,y_3).
\end{align*}
where $\epsilon=\pm1$, or
\item[{\rm (b)}] if $0$ is a triple point then $f_t$ is of the form
\begin{align*}
f_t(x)&=(x_1,x_2,x_3,0,0),\\
f_t(y)&=(y_1,t,0,y_2,y_3),\\
f_t(z)&=(0,z_2,z_1,z_2,z_3).
\end{align*}
\end{itemize}\qed
\end{prp}
If $\epsilon=+1$ ($\epsilon=-1$) in cases (a) in the above propositions
we say 
that $0$ is an {\em elliptic self tangency point} 
(a {\em hyperbolic self tangency point}) of $f_0$.  
\subsection{The stratum of codimension two}\label{vdf2}
\noindent
To study first order invariants of generic immersions (see Section
~\ref{1stord}) we need a 
description of the 
codimension two stratum of the discriminant hypersurface. We restrict
attention to immersions into 5-space.   

\begin{prp}\label{verdef52}
The codimension two stratum $\Si_5^2\subset\Si$ is the set of all
immersions $f_0\colon S^3\to\tR^5$ such that exactly one of the
following holds  
\begin{itemize}
\item[{\rm (a)}] $f_0$ has two distinct self tangency points,
\item[{\rm (b)}] $f_0$ has two distinct triple points,
\item[{\rm (c)}] $f_0$ has one self tangency point and one triple point,
\item[{\rm (d)}] $f_0$ has a degenerate self tangency point at 
$0=f(p)=f(q)$, in  which case its versal deformation $f_{t,s}$ is of
the form 
\begin{align*}
f_{t,s}(x)&=(x_1,x_2,0,x_3,0),\\
f_{t,s}(y)&=(y_1,y_2,y_1^2+y_2(y_2^2+s)+t,0,y_3),
\end{align*}
up to choice of coordinates $x$, $y$ around $p$, $q$, respectively,
and coordinates in $\tR^5$. 
\end{itemize}\qed
\end{prp}
The versal deformations of (a)-(d) are evident. They are just products
of 1-dimensional versal deformations.

\begin{figure}[htbp]
\begin{center}
\includegraphics[angle=0, width=10cm]{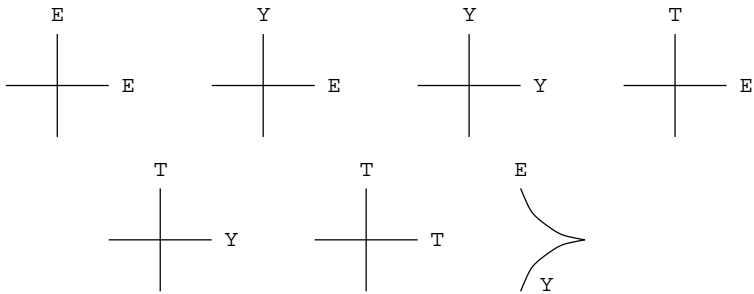}
\end{center}
\caption{The discriminant hypersurface intersected with a generic
2-disk}\label{figcod2} 
\end{figure}

Let $E$, $Y$ and $T$ denote the codimension one parts of $\Si_5$
consisting of all immersions with one
elliptic self tangency point, one hyperbolic self tangency point, and
one triple point, respectively. Figure ~\ref{figcod2}, which is a
consequence of Proposition ~\ref{verdef52}, shows the
possible intersections of the discriminant 
hypersurface $\Si_5$ and a small 2-disk in $\tF_5$ 
which meets $\Si_5^2$ and is transversal to $\Si_5$.
\subsection{Proofs}\label{discpf}
\noindent
The proofs of the propositions in Section ~\ref{vdf1} and Section ~\ref{vdf2}
are all similar: First we need to find coordinates close to the
preimages of the exceptional point such that the map is given by the
expression stated. Then we show that the deformation given is
infinitesimally versal. A standard theorem in singularity theory then
implies that the deformation is versal. As an example we write out
the proof of Proposition ~\ref{verdef41} (a): 

Since $f_0$ has a tangency at $0$ the tangent planes of the two sheets
$X$ and $Y$
meeting there must agree. We may assume that this tangent plane is the plane
of the first three coordinates. By the implicit function theorem we
can choose coordinates so that the map of the first sheet is
$$
f_0(x)=(x_1,x_2,x_3,0).
$$
We may now look upon the second sheet as the graph of a function
$\phi\colon\tR^3\to\tR$. The requirement that $Y$ is tangent to $X$ at
$0$ implies that $\phi(0)=0$ and $\frac{\partial\phi}{\partial
y_i}(0)=0$, $i=1,2,3$. Moreover, the second partials must be
nondegenerate since the immersion is in the codimension 1 part of the
discriminant hypersurface. Changing coordinates in $\tR^4$, by adding
a function of the first three coordinates that vanishes to order 3 to the
fourth, we may  
assume that $\phi$ is identical to the second order terms of its
Taylor polynomial. Choosing the coordinates in $Y$ appropriately then
gives $f_0$ the desired form.  

We must check infinitesimal versality. If $(z_1,\dots,z_4)$ are
coordinates on $\tR^4$ then this amounts to showing that 
any smooth variations
$$
\alpha(x)=\sum_{i_1}^4\alpha_i(x)\frac{\partial}{\partial z_i},\quad
\beta(x)=\sum_{i_1}^4\beta_i(y)\frac{\partial}{\partial z_i},
$$
can be written as
\begin{align*}
\alpha(x)&=\sum_{i=1}^3 a_i(x)\frac{\partial f_0}{\partial
x_i}+k(f_0(x))+c\dot f_t(x)|_{t=0},\\
\beta(y)&=\sum_{i=1}^3 b_i(y)\frac{\partial f_0}{\partial
y_i}+k(f_0(y))+c\dot f_t(y)|_{t=0},
\end{align*}
where $k$ is a vector field on $\tR^4$ and $c$ is a constant. Writing
these equations out we get the system
\begin{align*}
\alpha_i(x)&=a_i(x)+k_i(x,0),\quad i=1,2,3,\\
\alpha_4(x)&=k_4(x,0),\\
\beta_i(y)&=b_i(y)+k_i(y,y_1^2+y_2^2+\epsilon y_3^2),\quad i=1,2,3,\\
\beta_4(y)&=2(y_1b_1(y)+y_2b_2(y)+y_3b_3(y))+k_4(y,y_1^2+y_2^2+\epsilon
y_3^2)+c.
\end{align*}
Let $k_4(z)=\alpha_4(z_1,z_2,z_3)$. Then the fourth equation holds. Let
$c=\beta_4(0)-k_4(0)$. Let  $\phi(y)=\beta_4(y)-k_4(f_0(y))-c$ and
choose $c$ so that $\phi(0)=0$. If $y_i\ne 0$ let 
\begin{align*}
b_1(y)&=\frac{1}{6y_1}
\left(\phi(y_1,y_2,y_3)-\phi(0,y_2,y_3)+\phi(y_1,0,y_3)-
\phi(0,0,y_3)+\phi(y_1,0,0)\right),\\ 
b_2(y)&=\frac{1}{6y_2}
\left(\phi(y_1,y_2,y_3)-\phi(y_1,0,y_3)+\phi(y_1,y_2,0)-
\phi(y_1,0,0)+\phi(0,y_2,0)\right),\\ 
b_3(y)&=\frac{1}{6y_3}
\left(\phi(y_1,y_2,y_3)-\phi(y_1,y_2,0)+\phi(0,y_2,y_3)-
\phi(0,y_2,0)+\phi(0,0,y_3)\right),
\end{align*}
and if $y_i=0$ let
\begin{align*}
b_1(y)&=\frac{1}{6}
\left(\frac{\partial\phi}{\partial
y_1}(0,y_2,y_3)+\frac{\partial\phi}{\partial y_1}(0,0,y_3)+
\frac{\partial\phi}{\partial y_1}(0,0,0)\right),\\ 
b_2(y)&=\frac{1}{6}
\left(\frac{\partial\phi}{\partial
y_2}(y_1,0,y_3)+\frac{\partial\phi}{\partial y_2}(y_1,0,0)+
\frac{\partial\phi}{\partial y_2}(0,0,0)\right),\\ 
b_3(y)&=\frac{1}{6}
\left(\frac{\partial\phi}{\partial
y_3}(y_1,y_2,0)+\frac{\partial\phi}{\partial y_3}(0,y_2,0)+
\frac{\partial\phi}{\partial y_3}(0,0,0)\right).
\end{align*}
Then the last equation holds. Choosing first $k_i(z)=k_i(z_1,z_2,z_3)$,
$i=1,2,3$
so that the three remaining $\beta$-equations hold and then choosing
$a_i$, $i=1,2,3$ so that the remaining $\alpha$-equations hold, infinitesimal
versality is proved.

\section{Morse modifications, linking, and first order invariants}
\noindent
In this section,
we study how the self intersections of  generic immersions and their
preimages transform as the immersions cross the discriminant
hypersurface. These transformations give rise to invariants of 
immersions: 
A function on $\tF_n$, which is constant on path
components will be called an {\em invariant of regular homotopy}.      
A function on $\tF_n-\Si_n$ which is constant on path
components of $\tF_n-\Si_n$ is an {\em invariant of generic
immersions}. In our study of invariants of generic immersion we are
interested in how they change when we pass $\Si_n$. This is described
in terms of jumps:   

Let $f_t\colon S^3\to\tR^n$ be a path in $\tF_n$, $n=4,5$ intersecting
$\Si_n^1$ transversally at $f_0$. Let $\delta>0$ be small. Let $v$ be
an invariant of generic immersions. Then  
$$
\nabla v(f_0)=v(f_\delta)-v(f_{-\delta}),
$$ 
is a locally constant function on $\Si_n^1$, defined up to sign. We
call it {\em the jump of} $v$. In Section ~\ref{1stord}, we get rid of
the sign ambiguity in the definition of $\nabla v$. We prove
Theorem ~\ref{thmJSt} and end the section with an axiomatic
characterization of strangeness $\St$ for immersions $S^3\to\tR^5$.

\subsection{Morse modifications}\label{secmorinv}
\noindent
Throughout this section, let $f_t$ be a path in $\tF_n$ intersecting
$\Si_1$ transversally at $f_0$ and let $\delta\ge 0$ be small enough
so that $f_t$ is generic for $0<|t|\le\delta$.

Let $n=5$ and let $f_0$ have a self tangency point. By Proposition
~\ref{verdef51}, the self 
intersection of $f_\delta$ is obtained from the self intersection of
$f_{-\delta}$ by a single Morse modification. Define 
$$
J(f)=\text{the number of connected components of $M_f$},
$$ 
for generic immersions $f$. 
\begin{lma}\label{1stJ}
$J$ is an invariant of generic immersions $S^3\to\tR^5$. It
jumps by $\pm 1$ when crossing the self tangency part of $\Si_5^1$ and
remains constant when crossing the triple point part.
\end{lma}
\begin{pf}
Immediate from Proposition \ref{verdef51}.
\end{pf}

Let $n=4$ and let $f_0$ have an exceptional quadruple point. By
Proposition ~\ref{verdef41} the number of quadruple points of $f_\delta$
differs from the number of quadruple points of $f_{-\delta}$ by 
$\pm2$. Define 
$$
Q(f)=\text{the number of quadruple points of $f$},\text{ and }
Q_2(f)= Q(f)\mod{2}\in\tZ_2,
$$
for generic immersions $f$.  
\begin{prp}\label{mor4}
$Q$ is an invariant of generic immersions $S^3\to\tR^4$. It
jumps by $\pm2$ when crossing the part of
$\Si_4^1$ that consists of immersions with an exceptional quadruple point
and does not jump when crossing any other part of
$\Si_4^1$. Furthermore, $Q_2$ is an invariant of
regular homotopy.
\end{prp}
\begin{pf}
The first part is immediate from Proposition \ref{verdef41}. Let $f$ and $g$
be regularly homotopic then $f$ and $g$ can be joined by a path in
$\tF_4$ which intersects $\Si_4^1$ transversally. Since $Q$ jumps by
$\pm2$ or $0$ on such intersections, the lemma follows.
\end{pf}

Let $n=4$ and let $f_0$ have an exceptional triple point. 
By
Proposition ~\ref{verdef41}, $C_{f_\delta}$ is 
obtained from $C_{f_{-\delta}}$ by a single
Morse modification (for notation, see Section ~\ref{genim}). 
Define 
$$
T(f)=\text{the number of components of $C_f$},
$$
for generic immersions $f$. 
\begin{prp}
$T$ is an invariant of generic immersions $S^3\to\tR^4$. It
jumps by $\pm1$ when crossing the part of
$\Si_4^1$ that consists of immersions with an exceptional triple
point and does not jump when crossing any other part of 
$\Si_4^1$.
\end{prp}
\begin{pf}
Immediate from Proposition \ref{verdef41}.
\end{pf}

Let $n=4$ and let $f_0$ have an exceptional double point (a self
tangency point). 
By
Proposition ~\ref{verdef41}, $F_{f_\delta}$ is obtained from 
$F_{f_{-\delta}}$  by a single
Morse modification (for notation, see Section ~\ref{genim}). 
Define 
$$
D(f)=\chi(F_f)\quad\text{ and }\quad
D_2(f)= D(f)\mod{2}\in\tZ_2,
$$
for generic immersions $f$, where $\chi$ denotes the Euler characteristic. 
\begin{prp}
$D$ is an invariant of generic immersions $S^3\to\tR^4$. It
jumps by $\pm2$ when crossing the part of
$\Si_4^1$ that consists of immersions with an exceptional double point 
and does not jump when crossing any other part of
$\Si_4^1$. Furthermore, $D_2$ is an invariant of
regular homotopy.
\end{prp}
\begin{pf}
Similar to the proof of Proposition ~\ref{mor4}.
\end{pf}

\subsection{Linking}\label{seclk}
\noindent
We construct an invariant of generic immersions $S^3\to\tR^5$ which
jumps under triple point moves and stays constant under self tangency moves:
Let $f\colon S^3\to\tR^5$ be a generic immersion with self
intersection $M_f$. Consider the preimage ${\wt
M}_f=f^{-1}(M_f)$. 

Choose a normal vector field $w$ along ${\wt M}_f$ satisfying the
following condition: If
${\wt M}_f'$ denotes the result of pushing ${\wt M}_f$ slightly along
$w$ we require that 
\begin{equation*}
[{\wt M}_f']=0\in H_1(S^3-{\wt M}_f).
\tag{$\dagger$}
\end{equation*}
Any two vector fields satisfying this condition are homotopic. For
existence, note that we can take $w$ as the normal vector field of
the boundary in a Seifert-surface of the link ${\wt M}_f\subset S^3$.

We define a vector field $v$ along $M_f$: For $p\in M_f$, let 
$$
v(p)=df\,w(p_1)+df\,w(p_2),\text{ where }\{p_1,p_2\}=f^{-1}(p).
$$
Let $M_f'$ denote the result of pushing $M_f$ slightly along $v$. Then
$M_f'\subset\tR^5-f(S^3)$.

By Alexander duality 
$H_1(\tR^5-f(S^3))\cong H^3(f(S^3))$ and $H^3(f(S^3))\cong\tZ$ since
there is a triangulation of $S^3$ giving a
triangulation of $f(S^3)$ after identifications in the $0$ and $1$
skeletons only. An
orientation of $S^3$ gives a canonical generator of
$H_1(\tR^5-f(S^3))$: The boundary of a small 
2-disk intersecting $f(S^3)$ transversally in one point and
oriented in such a way that the intersection number is positive. 
Note that this intersection number is independent of the ordering of
the 2-disk and $f(S^3)$.

Recall that there is an induced orientation on $M_f$ (Section
~\ref{genim}) and hence on $M_f'$. Define
$$
\lk(f)=[M_f']\in H_1(\tR^5-f(S^3))=\tZ\quad\text{ and }\quad
\lambda(f)=\lk(f)\mod{3}\in\tZ_3.
$$
Note that $\lk(f)$ is well defined since an homotopy between two vector 
fields $w$ and $w'$ satisfying ($\dagger$)
induces a homotopy between the shifted self intersections in $\tR^5-f(S^3)$. 

\begin{lma}\label{lk}
$\lk$ is an invariant of generic immersions. It jumps by $\pm3$ when
crossing the triple point part of $\Si_5^1$ and remains constant when
crossing the
self tangency part. Furthermore, $\lambda$ is an invariant of regular
homotopy. 
\end{lma}
\begin{pf}
The second part follows from the first exactly as in Proposition
~\ref{mor4}. We prove the first part: 

Suppose that $f_t$ is a path in $\tF_5$ intersecting $\Si_5^1$
transversally 
at $f_0$ and let $\delta>0$ be small. We restrict attention to a small
neighborhood of the exceptional multiple point of $f_0$ since all the
immersions $f_t$, $|t|\le\delta$ can be assumed to agree outside of
this neighborhood. 

Assume that $f_0$ has an elliptic self tangency point $f_0(p)=f_0(q)$. Using
Proposition ~\ref{verdef51} we can write 
\begin{align*}
f_t(x) &=(x_1,x_2,t,x_3,0).\\
f_t(y) &=(y_1,y_2,y_1^2+ y_2^2,0,y_3),
\end{align*}
where $x$ are coordinates around $p$ and $y$ are coordinates around
$q$. The preimages of the newborn self intersection circle $c$ of
$f_\delta$ are $\{x_1^2+x_2^2=\delta, x_3=0\}$ and
$\{y_1^2+y_2^2=\delta, y_3=0\}$. We may assume that the field $w$ is
given by $w(x)=\frac{\partial}{\partial x_3}$ and
$w(y)=\frac{\partial}{\partial y_3}$. Now shift $c$ a small distance
$\epsilon$ along $df_\delta w(x)+df_\delta w(y)$ to get $c'$. Then $c'$
is the boundary of the 2-disk
$$
(r\delta\cos(v),r\delta\sin(v),\delta,\epsilon,\epsilon),\quad
0\le v\le2\pi,\, 0\le r\le 1,
$$ 
which does not intersect $f_\delta(S^3)$. Hence,
$\lk(f_\delta)=\lk(f_{-\delta})$. 

Assume that $f_0$ has a hyperbolic self tangency point. As above we
can write
\begin{align*}
f_t(x) &=(x_1,x_2,t,x_3,0).\\
f_t(y) &=(y_1,y_2,y_1^2- y_2^2,0,y_3).
\end{align*}
The preimages of the self intersection are $\{x_1^2-x_2^2=\delta,
x_3=0\}$ and $\{y_1^2-y_2^2=\delta, y_3=0\}$. We can choose the field
$w$ so that $w(x)=\frac{\partial}{\partial x_3}$ and
$w(y)=\frac{\partial}{\partial y_3}$ close to $p$ and $q$,
respectively. 
We
note that with vector fields as above and equal outside neighborhoods
of $p$ and $q$ the condition ($\dagger$) holds for $f_{-\delta}$ if 
and only if it holds for $f_\delta$. 
The rest of the argument is similar to the elliptic case. 

Assume that $f_0$ has a triple point.
According to Proposition ~\ref{verdef51} we can find coordinates $x$, $y$ and
$z$ centered at $p_1$, $p_2$ and $p_3$, respectively such that 
close to the triple point $f_0(p_1)=f_0(p_2)=f_0(p_3)$ we have
\begin{align*}
f_t(x)&=(x_1,x_2,x_3,0,0),\\
f_t(y)&=(y_1,t,0,y_2,y_3),\\
f_t(z)&=(0,z_2,z_1,z_2,z_3).
\end{align*}
Denote the neighborhoods of $p_1$, $p_2$ and
$p_3$ by $X$, $Y$ and $Z$, respectively. If we orient these by
declaring the frames
$(\partial_1,\partial_2,\partial_3)$ to have the positive orientation
then the oriented self intersections of $f_0$ are the lines 
\begin{align*}
f_t(X)\cap f_t(Y) &= (s,t,0,0,0),& s\in\tR,\\
f_t(X)\cap f_t(Z) &= (0,0,s,0,0),& s\in\tR,\\
f_t(Z)\cap f_t(Y) &= (0,t,0,t,-s),& s\in\tR.
\end{align*} 
The inverse image of the self intersection of $f_t$, $t<0$ is shown in
Figure ~\ref{figinvtrip}.
\begin{figure}[htbp]
\begin{center}
\includegraphics[angle=-90, width=10cm]{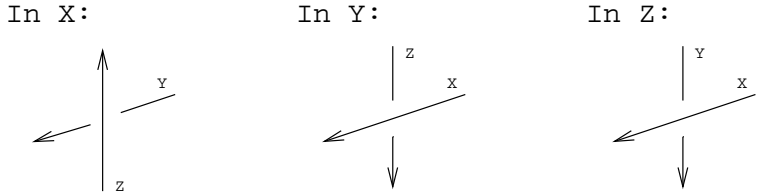}
\end{center}
\caption{Crossings close to a triple point}\label{figinvtrip}
\end{figure}

We calculate $\lk(f_\delta)$ in terms of $\lk(f_{-\delta})$ and
have to take two things into consideration:
First the motion of $Y$ itself. Second, the changes in the field $w$
that this motion causes. 

Assume that the shifting distance is
very small in comparison to $\delta$.
Let $l_t(X,Z)$ denote the intersection $f_t(X)\cap f_t(Z)$ shifted
along $v=df\,w(x)+df\, w(z)$. Let
$D_t(X,Z)$ denote a part of a disc in $\tR^5$ bounded by $l_t(X,Z)$. We may
assume that $D_t(X,Z)$ is a shift along an extension of $v$ of a disk
in $f_t(X)$.  

First, the
intersection $D_t(X,Z)\cap f_t(Y)$ does not depend on the field $w$
since the shift is assumed to be very small in comparison to $t$. We
calculate the 
change in $D_t(X,Z)\cap f_t(Y)$ from the preimage in $X$. From Figure
~\ref{figbiglk} it follows that  the algebraic number of intersection
points in $D_{-\delta}(X,Z)\cap f_{-\delta}(Y)$ differs from that
corresponding to $D_{\delta}(X,Z)\cap f_{\delta}(Y)$ by $+1$. 
\begin{figure}[htbp]
\begin{center}
\includegraphics[angle=-90,width=10cm]{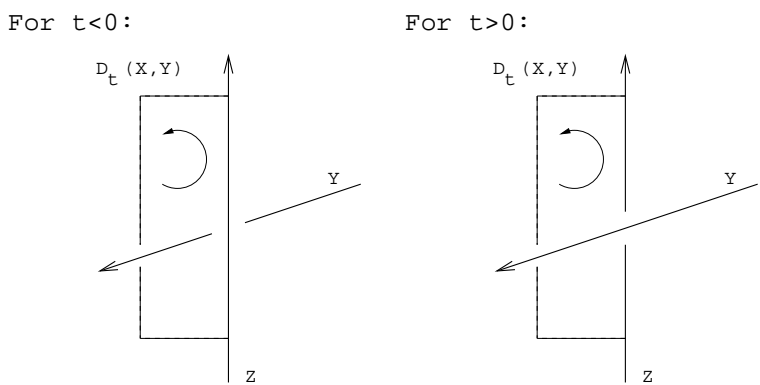}
\end{center}
\caption{$D_t(X,Y)$ intersected with $Y$.}\label{figbiglk}
\end{figure}

We now take the field $w$ into consideration. 
In Figure ~\ref{figseif} we see a Seifert-surface in $S^3$ of ${\wt
M}_f$ close to a crossing point.
\begin{figure}[htbp]
\begin{center}
\includegraphics[angle=-90,width=3cm]{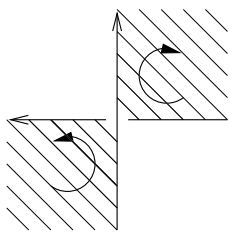}
\end{center}
\caption{A Seifert-surface.}\label{figseif}
\end{figure}
 
In Figure ~\ref{figsmalllk} 
we see $w$ chosen as the inward normal in a Seifert-surface of ${\wt
M}_{f_t}$. 
\begin{figure}[htbp]
\begin{center}
\includegraphics[angle=-90,width=10cm]{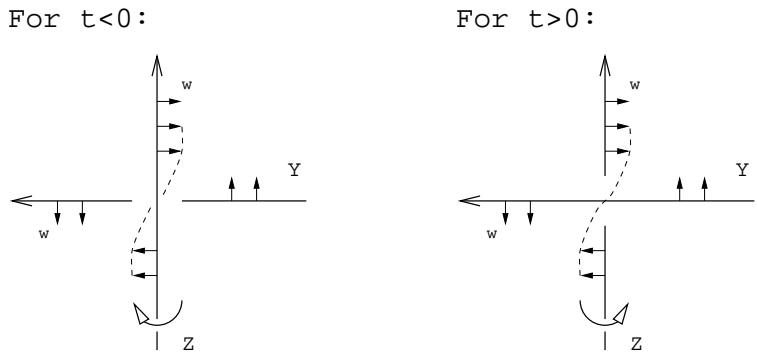}
\end{center}
\caption{Change of $w$ on passing a triple point.}\label{figsmalllk}
\end{figure}
As we move 
through the triple point the direction of rotation of the vector field $w$
is changed. This change in rotation gives rise to a new positive
intersection point in 
$D_t(X,Z)\cap f_t(Z)$ for $t>0$. As seen in Figure ~\ref{figtwlk}. 
\begin{figure}[htbp]
\begin{center}
\includegraphics[angle=-90,width=8cm]{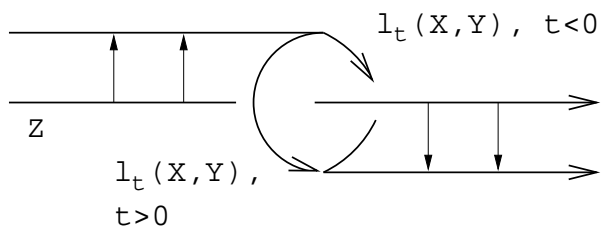}
\end{center}
\caption{$D_t(X,Z)$ intersected with $f_t(Z)$.}\label{figtwlk}
\end{figure}

Similarly, we
get a new positive intersection point in $D_t(X,Z)\cap f_t(X)$. Thus,
in total the algebraic number of intersection points in 
$D_t(X,Z)\cap (f_t(X)\cup f_t(Y)\cup f_t(Z))$ increases by   
$1$ when $t$ is changed from $-\delta$ to $\delta$.

Similarly, the intersection numbers corresponding to 
$D_t(X,Y)\cap (f_t(X)\cup f_t(Y)\cup f_t(Z))$ and 
$D_t(Y,Z)\cap (f_t(X)\cup f_t(Y)\cup f_t(Z))$ increase by
$1$ when $t$ is changed from $-\delta$ to $\delta$. 

It follows that $\lk(f_\delta)-\lk(f_{-\delta})=3$.
This proves the lemma.
\end{pf}

\begin{cor}\label{homlambda}
For $\xi$ in $\Imm$, let $\lambda(\xi)=\lambda(f)$, where $f$ is a
generic immersion representing $\xi$. Then
$$
\lambda\colon\Imm\to\tZ_3,
$$
is a homomorphism.
\end{cor}
\begin{pf}
As a function, $\lambda$ is well defined by Lemma ~\ref{lk}. Clearly,
it is additive under connected sum and thus a homomorphism.
\end{pf}
\begin{rmk}\label{3signs}
In knot theory one assigns a sign to crossings as in Figure
~\ref{figinvtrip} by comparing 
the orientation given by $(v_1,v_2,v_{12})$ to that of the ambient
space, where $v_1$ is the tangent vector to the first branch, $v_2$
the tangent vector of the second, and $v_{12}$ a vector from the second
to the first branch.
We note that all crossings appearing in Figure ~\ref{figinvtrip} are
positive. Changing the 
orientation of one of the sheets changes the sign of all
crossings. Hence, {\em all three crossings appearing close to a triple
point always have the same sign.} 
\end{rmk}

\subsection{Finite order invariants}\label{1stord}
\noindent
In this section we summarize the properties of finite order invariants
of generic  immersions that are needed to prove Theorem ~\ref{thmJSt}.
Recall that we have defined the jump $\nabla v$ of an invariant $v$ of
generic immersions. It was defined only up to sign. To get rid of this
sign we need a coorientation of $\Si_n^1$: 

If $f_0\in\Si_n^1$ then
there is a neighborhood $U$ of $f_0$ in $\tF_n$ which is cut into two
parts by $\Si_n^1$. We make a choice of a positive part and a negative
part of $U$. A coherent choice like that for all $f_0\in\Si_n^1$ is a
{\em coorientation of $\Si_n$}. A coorientation enables us to make
$\nabla v$ well defined: In the definition
$$
\nabla v(f_0)=v(f_\delta)-v(f_\delta),
$$
we require that $f_\delta$ is on the positive and $f_{-\delta}$ on the
negative side of $\Si_n^1$ at $f_0$. As mentioned, $\nabla v$ is then
locally constant on $\Si_n^1$ and we may consider $\nabla v$ as an
element of $H^0(\Si_n^1)$. 

Let $f_t$ be a path in $\tF_5$ intersecting $\Si_5^1$ transversally at
$f_0$. We coorient $\Si_5$ as follows: 
If $f_0\in E\cup Y$ then $f_\delta$ is on the positive side of $E\cup
Y$ at $f_0$ if $M_{f_\delta}$ has more components than
$M_{f_{-\delta}}$. If $f_0\in T$ then $f_\delta$ is on the positive
side of $T$ if $\lk(f_\delta)\ge\lk(f_{-\delta})$. 

It is easy to check that this coorientation is {\em continuous} (see
\cite{E2}, Section 7.6). This enables us to iterate the above
construction and define inductively 
$$
\nabla^{k+1}v(f_{0,0,\dots,0})=
\nabla^k(f_{\delta,0,\dots,0})-\nabla^k(f_{-\delta,0,\dots,0}).
$$
Here $f_{0,\dots,0}$ is an immersion with $k+1$ distinct degenerate
points (self tangencies or triple points) 
and $f_{t,0,\dots,0}$ is a path in $\Si_5^k\cup\Si_5^{k+1}$, where
$\Si_5^j$ is the space of
immersions with $j$ distinct degenerate points, intersecting
$\Si_5^{k+1}$ transversally. Then $\nabla^k v$ is an element in
$H^0(\Si_5^k)$ (see \cite{E1}, Remark 9.1.2). 
We say that {\em an invariant $v$ is of order $k$} if $\nabla^{k+1}v\equiv
0$.

Finally, we remark that the space of invariants of generic immersions,
i.e.\linebreak 
$H^0(\tF_n-\Si_n)$, splits as a direct sum over the path components
of $\tF_n$. That is $H^0(\tF_n-\Si_n)=\bigoplus_{U}H^0(U-(\Si_n\cap U))$,
where $U$ runs over the path components of $\tF_n$.

\subsection{Proof of Theorem ~\ref{thmJSt}}\label{pfJSt}
\noindent
By Lemmas ~\ref{1stJ} and ~\ref{lk}, $J$ and $L$ are
invariants which changes as claimed. Clearly, $L$ and $J$ are independent. 

Since $\nabla J$ ($\nabla L$)
is $1$ ($0$) on all self tangency parts of $\Si_5^1$ and is $0$
($1$) on  all triple point parts of $\Si_5^1$, it follows from 
Proposition ~\ref{verdef52} (see Figure ~\ref{figcod2}) that $\nabla^2
J\equiv 0$ and $\nabla^2 L\equiv 0$. Thus, they both have order one. 

Let $U$ be a path component of $\tF_5$ and $v$ an invariant, the jump
of which is constant on $T\cap U$ and $(E\cup Y)\cap U$ (for notation
see Section ~\ref{vdf1}). Then the 
jump of $v|U$ 
is a linear combination of the jumps of $J|U$ and $L|U$ and it follows
that $v|U$ is (up to a constant term) a linear combination of $J|U$ and
$L|U$. Thus, the theorem follows once we show that any first order
invariant has this property.

Let $f$ and $g$ be regularly homotopic immersions with one self
tangency point of the same 
kind each or with one triple point each. Using a diffeotopy of $\tR^5$
we can move the exceptional point of $g$ to that of $f$ and using a
diffeotopy of $S^3$ we can make the maps agree in small disks around
the preimages of this point. Choose unknotted paths in $S^3$
connecting these disks. Since there are no 1-knots in $\tR^5$, we can
after diffeotopy of $\tR^5$ assume that the maps agree also in a
neighborhood of these arcs. After this is done we have two immersions
that agree on a 3-disk in $S^3$ and are regularly homotopic. We must
show that they are regularly homotopic keeping this 3-disk
fixed. However, by the Smale-Hirsch h-principle the obstruction to this
regular homotopy is the difference of the Smale invariants which is
zero. Choosing a generic regular homotopy we see that we can join $f$
to $g$ with a path in $\Si^1_5\cup\Si_5^2$ intersecting $\Si_5^2$
transversally only at immersions with two distinct exceptional points.

Let $v$ be a first order invariant. Then $\nabla^2v\equiv0$. This and
the fact that $\nabla v$ is locally constant on $\Si_5^1$ implies that
$\nabla v$ remains constant along a path as above. Hence,
if $f$ and $g$ both belong to either one of $E\cap U$, $Y\cap U$, or
$T\cap U$ then $\nabla v(f)=\nabla v(g)$. 

It remains to show that $\nabla v|E\cap U$ equals $\nabla v|Y\cap U$:
Consider the last part of Figure ~\ref{figcod2} in Section
~\ref{vdf2}. Let $f$ be an 
immersion in the $E$-branch and $g$ an immersion in the $Y$-branch. Pick
a small loop $h_t$, $0\le t \le 1$, $h_1=h_0$ around the cusp such
that it intersects the $E$-branch in $f$ and the $Y$-branch in $g$. Then
$0=v(h_1)-v(h_1)=\nabla v(f)-\nabla v(g)$. Hence, 
$$
\nabla v|E\cap U\equiv\nabla v(f)=\nabla v(g)\equiv\nabla v|Y\cap U.
$$
The theorem follows. \qed

\subsection{$L$ and $\St$}\label{secLSt}
\noindent
Arnold's invariant $\St$ for plane curves
(see \cite{A2}) can be characterized axiomatically as an invariant
which jumps by $1$ on triple points, does not jump under self
tangencies, and is additive under connected summation. 

The invariant $L=\frac13(\lk+{\wt\lambda})$ jumps by one on triple
points but is not additive under connected sum. The reason for this is
that the invariant ${\wt \lambda}$ of order zero is not additive under
connected sum.

The invariant $\lk$ is additive under connected sum. It jumps by
$\pm3$ on $T$ and not at all on $E\cup Y$. Thus, $\lk$ and $J$
generates the space of $\tQ$-valued invariants. 

We want an invariant which is additive under connected sum and
together with $J$  generates the space of $\tZ$-valued invariants.
To accomplish this we note that, either 
$\lk+\Omega$ or $\lk-\Omega$ is divisible by $3$, where $\Omega$ is the
Smale invariant, which we may consider integer valued after the choice 
of generator of $\pi_3(V_{5,3})\cong\tZ$ made in Lemma ~\ref{pi3V53}.

In Lemma ~\ref{lkgen} we calculate $\lk(f)$ for 
an immersion with $\Omega(f)=1$ and the result is $\lk(f)=2$. 
Since $\lambda=\lk\mod{3}$ is an invariant of regular homotopy it
follows that
$\frac13(\lk(f)+\Omega(f))$ is an integer valued invariant which jumps
by $\pm 1$ on triple points, does not jump on self tangencies, and is
additive under connected sum (since both $\lk$ and $\Omega$ are). We
therefore define 
$$
\St(f)=\frac13((\lk(f)+\Omega(f))).
$$

To complete the analogy with Arnold's $\St$, we need to define the
coorientation of the triple point stratum $T$ of $\Si_5^1$ in terms of the
local picture near a triple point. Recall that all crossings in the
preimage close to a triple point have the same sign (Remark
~\ref{3signs}). By Lemma ~\ref{lk} our coorientation of $T$ agrees
with the following one defined in local terms:

{\em We say that we pass $T$ in the positive direction if the crossing
signs in the preimage close to the triple point changes from positive
to negative.} 

\begin{prp}\label{Strange}
There is a unique invariant $\St$ such 
that it jumps by $+1$ ($-1$) under
positive (negative) triple point moves, it does not jump under self tangency
moves, it is additive under connected sum, its value on an immersion
changes sign if the immersion is composed with an 
orientation reversing diffeomorphism of $S^3$, and takes the values $1$ on
$f_{-\frac12}$ (see Proposition ~\ref{theg}).
\end{prp}
\begin{pf}
A generic immersion $g$ with $\Omega(g)\ne 0$ is regularly homotopic
to a connected sum of  
$|\Omega(g)|\in\tZ_+$ copies of $f_{-\frac12}$ or $f_{-\frac12}\circ
r$, where $r$ is an orientation reversing diffeomorphism.
A generic immersion $g$ with $\Omega(g)=0$ is regularly homotopic to 
$f_{\frac12}\star f_{-\frac12}$. 

We know the values of our invariants on
connected sums of $f_{\pm\frac12}$. Thus, we know the value on any
other generic immersion by adding the jumps along a path connecting it
to such a 
connected sum. This proves uniqueness. 

For existence, we need to show
that $\St$ changes sign if the orientation is reversed. For $\Omega$
we already know this (Lemma ~\ref{invrs}). To see that the same is
true for $\lk$ note that the induced orientation of $M_f$ does not
change if we reverse the orientation on $S^3$. The orientation of
$f(S^3)$ does change and hence $\lk(f)=[M_f']\in H_1(\tR^5-f(S^3))$
does change sign.
\end{pf}

\section{Pin structures and twist framings}
\noindent
In this section we show that there is an induced pin structure on the
self intersection surface of a generic immersion $S^3\to\tR^4$. The Brown
invariant of this pin structure is unchanged under regular
homotopy. Actually an even stronger result holds, if two
generic immersions into $\tR^4$ are regularly homotopic in $\tR^5$
after composing them with the inclusion 
$\tR^4\to\tR^5$ then the corresponding Brown invariants are equal. 

Using the geometry of self intersections of generic immersions
$S^3\to\tR^5$, we define a $\tZ_4$-valued invariant of regular
homotopy. 

\subsection{Immersions into $\tR^4\subset\tR^5$}
\noindent
Let $g\colon S^3\to\tR^4$ be an immersion.
Its Smale invariant $\Omega(g)$ is an element of
$\pi_3(SO(4))=\tZ[\sigma]\oplus\tZ[\rho]$ (see Lemma
~\ref{pi3V43}). The normal bundle of $g$ is 
1-dimensional and orientable. Hence, it is trivial and $g$ admits a
normal field $n\colon S^3\to S^3$, such that the frame $(n(x),dfQ(x))$ (see
Section ~\ref{calc}) gives the positive orientation of $\tR^4$. If
$\Omega(g)=m[\sigma]+n[\rho]$ 
then, from Lemma ~\ref{pi3V43}, it follows that $m+1=\deg(n)$, the normal
degree of $g$. 

Let $s\colon S^3\to\tR^4$ be
the standard embedding and $\widehat s=s\circ r$, where $r$ is a
reflection in a hyperplane as in Lemma ~\ref{invrs}.
\begin{lma}
$$
\Omega({\widehat s})=-2[\sigma]+[\rho]\in\pi_3(V_{4,3})=\pi_3(SO(4)).
$$
\end{lma}
\begin{pf}
Let $i\colon\tR^4\to\tR^5$ be the inclusion. Then $i\circ s$ and
$i\circ{\widehat s}$ are regularly homotopic in $\tR^5$.
Hence, by Lemma ~\ref{Smin4in5}, $\Omega({\widehat s})\in N\subset SO(4)$. 

The normal degree of $\widehat s$ 
equals the normal degree 
of $s$ with opposite sign. Thus, $\Omega_4({\widehat
s})=-2[\sigma]+n[\rho]$. But then, $\Omega_4({\widehat s})\in N\subset
SO(4)$ implies that $n=1$. 
\end{pf}
Let $g,h\colon S^3\to\tR^4$ be immersions. 
\begin{prp}\label{4ss5}
If $i\circ g$ and $i\circ h$ are regularly homotopic in $\tR^5$ then either
$g$ and $h\star{\widehat s}\star\dots\star{\widehat s}$ or
$h$ and $g\star{\widehat s}\star\dots\star{\widehat s}$ are regularly
homotopic in $\tR^4$.
\end{prp}
\begin{pf}
If $i\circ g$ and $i\circ h$ are regularly homotopic in $\tR^5$ then
$\Omega(g)-\Omega(h)\in N$. Now, $\Omega({\widehat s})$ is a
generator of $N$ and the Smale invariant is additive under connected
sum. 
\end{pf}

\subsection{Functionals on curves in self intersections}\label{secfun}
\noindent
The invariants we are about to define originates from functionals on
curves in self intersections.
Before we can construct these functionals we need some preliminaries: 

Let $\tI$ denote the identity matrix.
\begin{dfn}
Define $\tA(4)\in SO(4)$ and $\tA(5)\in SO(5)$ by
$$
\tA(4)=\begin{bmatrix}
1&0&0&0\\
0&-1&0&0\\
0&0&0&1\\
0&0&1&0\\    
\end{bmatrix},\quad
\tA(5)=\begin{bmatrix}
1&0&0&0&0\\
0&0&0&1&0\\
0&0&0&0&1\\
0&1&0&0&0\\    
0&0&1&0&0\\    
\end{bmatrix}.
$$
\end{dfn} 
The group $\{\tI,\tA(n)\}$ acts on $SO(n)$ from the right. Denote the
corresponding quotient space $SO(n)/\tA(n)$ and let $p$ be the
projection.
\begin{prp}\label{pi_1(SO/A)}
For $n=4$ and $n=5$,
$$
\pi_1(SO(n)/\tA(n))=\tZ_4. 
$$
\end{prp}
\begin{pf}
See \cite{E1}, Proposition 5.2.
\end{pf}

Let $u\in\pi_1(SO(n))$ denote the nontrivial element.

The functional on curves in a self intersection surface of a generic
immersion into 4-space is constructed
as follows (for notation, see Proposition ~\ref{FgCg}):  

Let $g\colon S^3\to\tR^{4}$ be a generic immersion. Let $c$ be a
closed curve in $F_g$ meeting $F_g^0\cup F_g^1$ transversally
and let $c_S=p^{-1}(c)$. Note that $c_S$ is either a
union of two disjoint circles or one circle. 

Choose a
parameterization $r\colon I\to c$ of $c$ and a normal vector $n\colon
I\to TF_g$ of $c$ in $TF_g$. Then $r$ lifts to two parameterizations $r_1$ and 
$r_2$. If $c_S$ is connected then the path product $r_1\ast r_2$ is a
parameterization of $c_S$. 

Choose normal fields $\nu_1$ and $\nu_2$
of $\wM_g$ along $s\circ r_1$ and $s\circ r_2$, respectively. 
Let $n_1$ and $n_2$ be such that $dp(n_i)=n$. Assume that $\nu_i$ is
chosen so that $(ds\,\dot r_i,ds\,n_i,\nu_i)$ is a positively oriented
frame, $i=1,2$.  
Then $(dt\,\dot r,dt\,n,dg\,\nu_1,dg\,\nu_2)$ is a framing or twist
framing $X$ of $T\tR^{4}$ along $t\circ c$. 

Let $\vec c_S$ denote the loop
(or loops) in $SO(TS^3)$ (the principal $SO(3)$-bundle of the tangent
bundle of $S^3$) represented by $s\circ c_S$ with the framing $(ds\,\dot
r_i,ds\,n_i,\nu_i)$ along $s\circ r_i$, $i=1,2$ and let $[\vec c_S]$
denote the corresponding element in $H_1(SO(TS^3);\tZ_2)$. 

Let $[c,r,X]\in\pi_1(SO(4)/\tA(4))$
denote the homotopy class of the loop induced by $t\circ c$ with 
(twist) framing $X$ as above (which we may have to
orthonormalize). Then define   
\begin{dfn}\label{om_g}
$$
\omega_g(c)=\langle\xi_S,[\vec c_S]\rangle\cdot
p_\ast(u)+[c,r,X]+p_\ast(u)\in \pi_1(SO(4)/\tA(4)),
$$
where $\xi_S\in H^1(SO(TS^3);\tZ_2)$ denotes the unique spin structure on
$S^3$.
\end{dfn}
\noindent
This definition is independent of all choices (see \cite{E1}, Lemma 7.6). 

We extend the functional $\omega_g$ to collections
$L=\{c_1,\dots,c_m\}$ of oriented 
closed curves in $F_g$ transversal to $F_g^0\cup F^1_g$:
$$
\omega_g(L)=\sum_{i=1}^m\omega_g(c_i).
$$

Let $L$ be a collection of oriented closed curves in $F_g$ with $m$ transverse
intersections and $L'$ be the collection obtained from smoothing each
intersection, respecting orientation. Then
\begin{equation*}
\omega_g(L')=\omega_g(L)+m\cdot p_\ast(u),
\tag{$\dagger$}
\end{equation*}
and if $K$ and $L$ are isotopic collections then 
\begin{equation*}
\omega_g(K)=\omega_g(L),
\tag{$\ddagger$}
\end{equation*}
see \cite{E1}, Lemma 7.15.

The functional on self intersection circles of generic immersions
into $5$-space is constructed as follows:

Let $f\colon S^3\to\tR^5$ be a generic immersion. Let $c\subset M_f$
be a component and $c_S=f^{-1}(c)$. Note
that $c_S$ is either a union of two disjoint circles or one circle. 

Choose a
parameterization $r$ of $c$, agreeing with the orientation induced on
$c$ if $k$ is odd. Then $r$ lifts to two parameterizations $r_1$ and
$r_2$. If $c_S$ is connected then the path product $r_1\ast r_2$ is a
parameterization of $c_S$. 

Choose orthonormal framings $Y_1$ and $Y_2$ of
$N(c_S\subset S^3)$ (the normal bundle of $c_S$ in $S^3$) along $r_1$
and $r_2$ respectively. Then $(\dot r,dfY_1,dfY_2)$ is a framing or twist
framing $X$ of $T\tR^{5}$ along $c$. 

Let $\vec c_S$ denote the loop
(or loops) in $SO(TS^5)$ (the principal $SO(5)$-bundle of the tangent
bundle of $S^5$) represented by $c_S$ with the framing $(\dot
r_i,Y_i)$ and $[\vec c_S]$ the corresponding element in
$H_1(SO(TS^5);\tZ_2)$. Let $[c,r,X]\in\pi_1(SO(5)/\tA(5))$ denote the
homotopy class of the loop induced by $c$ with the (twist) framing
$X$ (which we may have to orthonormalize). Then define
\begin{dfn}
$$
\omega_f(c)=\langle\xi,[\vec c_S]\rangle\cdot
p_\ast(u)+[c,r,X]+p_\ast(u)\in \pi_1(SO(5)/\tA(5)),
$$
where $\xi\in H^1(SO(TS^5);\tZ_2)$ denotes the unique spin structure on
$S^5$.
\end{dfn}
\noindent
This definition is independent of all choices (see \cite{E1}, Lemma 6.5). 
\subsection{Pin structures on surfaces}
\noindent
Let $V$ be a vector space over $\tZ_2$ with a nonsingular symmetric
bilinear form $(x,y)\mapsto x\cdot y$. A $\tZ_4$-quadratic function
on $V$ is a function $q\colon V\to\tZ_4$ such that
$q(x+y)=q(x)+q(y)+2(x\cdot y)$, for all $x,y\in V$. 

There are four indecomposable $\tZ_4$-quadratic spaces:
\begin{alignat*}{3}
&\tP_+=(\tZ_2(a),\cdot,q),\quad & &a\cdot a=1,\quad & &q(a)=1,\\
&\tP_-=(\tZ_2(a),\cdot,q),\quad & &a\cdot a=1,\quad & &q(a)= -1,\\
&\tT_0=(\tZ_2(b)\oplus\tZ_2(c),\cdot,q),\quad & &b\cdot b=c\cdot c=0,\, b\cdot
c=1,\quad & &q(b)=q(c)=0,\\
&\tT_4=(\tZ_2(b)\oplus\tZ_2(c),\cdot,q),\quad & &b\cdot b=c\cdot c=0,\, b\cdot
c=1,\quad & &q(b)=q(c)=2.
\end{alignat*}

A $\tZ_4$-quadratic space $V$ is called {\em split} if it contains a subspace
$H$ such that $\dim H=\frac12\dim V$, $q|H=0$ and $H\cdot H=\{0\}$. Two
$\tZ_4$-quadratic spaces $V$ and $W$ belong to the same {\em Witt
class} if there exists spilt spaces $S_1$ and $S_2$ such that 
$V\oplus S_1\cong W\oplus S_2$. The Witt classes forms the {\em Witt group}
which is generated by $[\tP_+]$ (the Witt class of $\tP_+$) with
relations $8[\tP_+]=0$, $4[\tP_+]=[\tT_4]$, and $[\tP_+]+[\tP_-]=0$. 

Given a
$\tZ_4$-quadratic space $V$ with quadratic function $q$ we define
$$
\lambda(V,q)=\sum_{x\in V}e^{\frac{\pi i}{2}q(x)}.
$$
Then 
$$
\lambda(V,q)=\sqrt{2}^{\,\dim V}\left( \frac{1+i}{\sqrt{2}} \right)^m.
$$ 
Since $\frac{1+i}{\sqrt{2}}$ is an $8^{th}$ root of unity, $m$ modulo 8
is well defined. This is {\em Brown's invariant}. It is denoted $\beta(V,q)$
and gives an isomorphism between the Witt group and $\tZ_8$. We shall
sometimes write $\beta(q)$, dropping $V$. More details about
$\tZ_4$-quadratic spaces can be found 
in \cite{M}.

Let $F$ be a surface. By a pin structure on $F$ we shall mean
a $Pin^-$-structure on the tangent bundle $TF$ of $F$. 
There is a 1-1
correspondence between pin structures on $F$ and $\tZ_4$-quadratic
functions on $H_1(F;\tZ_2)$, see \cite{KT}, Theorem 3.2. 

\subsection{Invariance of the Witt class and a homomorphism
$\Imm\to\tZ_8$}\label{secbeta} 
\noindent
Let $g\colon S^3\to\tR^4$ be a generic immersion.
Choose, once and for all an isomorphism 
$$
\phi\colon\pi_1(SO(4)/\tA(4))\to\tZ_4.
$$
(This choice is discussed in Section ~\ref{secrelinv}.)
By equations ($\dagger$) and ($\ddagger$) in Section ~\ref{secfun} and
Lemma 3.4 in \cite{KT},
$\phi\circ\omega_g$ induces
a $\tZ_4$-quadratic function $q_g\colon H_1(F_g;\tZ_2)\to\tZ_4$ (and
hence a pin structure on $TF_g$). We define 
$$
\beta(g)=\beta(q_g),
$$ 
the Brown invariant of the $\tZ_4$-quadratic function $q_g$.  

\begin{lma}\label{breg4}
Let $g_0,g_1\colon S^3\to\tR^4$ be regularly homotopic generic
immersions. Then $\beta(g_0)=\beta(g_1)$
\end{lma}
\begin{pf}
Connect $g_0$ to $g_1$ by a path $g_t$ in $\tF_4$ which is transversal to
$\Si_4$. This means that $g_t$ intersects $\Si_4$ only transversally at
finitely many points in $\Si_4^1$. Clearly, $\beta$ remains constant
as long as we do not cross $\Si_4$ and we must show that $\beta$
remains unchanged when we cross $\Si^1_4$:

Passing $\Si_4^1$ at $g_{t_0}$,
where $g_{t_0}$ has a $j$-tuple point $j\ge 2$ does not
change $F_{g_t}$. It changes the preimage in $F_g$ of
$M_g^{4-j}$. Since the curves representing a basis of
$H_1(F_{g_t};\tZ_2)$, used to calculate $\beta$ can be chosen so that
they do not meet the discs where these changes occur, $\beta$ does not
change at such crossings. 

Passing $\Si_4^1$ at $g_{t_0}$, where $g_{t_0}$ has a degenerate
double point does change $F_{g_t}$ by a Morse modification. If it has
index $2$ or $0$ then $F_{g_t}$ is changed by addition or subtraction of an
$S^2$-component. This does not affect $\beta$. 
If it has index $1$ then $F_{g_t}$ is changed by addition or
subtraction of a handle. In this handle there is a newborn circle $c$,
$q_f(c)=0$ and $\beta$ remains unchanged, see \cite{E3}, Lemmas 5.2.1-8.
\end{pf}
Let $f\colon S^3\to\tR^5$ be an immersion. The normal bundle of $f$ is
trivial (Section ~\ref{calc}) which implies (\cite{Hi}, Theorem 6.4)
that $f$ is 
regularly homotopic to an immersion $f_1\colon S^3\to\tR^4$ composed
with the inclusion $\tR^4\to\tR^5$. Moreover, we may assume that $f_1$
is generic. Define 
$$
\beta(f)=\beta(f_1).
$$
\begin{prp}\label{hombeta}
For $\xi\in\Imm$, let $\beta(\xi)=\beta(f)$, where $f$ is an immersion
representing  $\xi$. Then
$$
\beta\colon\Imm\to\tZ_8,
$$
is a homomorphism.
\end{prp}
\begin{pf}
We show first that $\beta$ is well defined: 
If $f_1$ and $g_1$ are immersions into $\tR^4$ which are regularly
homotopic to $f$ in $\tR^5$ then $f_1$ is regularly homotopic to
$g_1\star{\widehat s}\star\dots\star{\widehat s}$ (Proposition
~\ref{4ss5}). We can perform the 
connected sum of two generic immersions into  $\tR^4$ in such a way
that the self intersection of the immersion obtained is the union of
the self intersections of the summands and at most two new $S^2$ self
intersection components. Hence, performing connected sum with the
embedding ${\widehat s}$ does not change $\beta$. Thus,
$\beta\colon\Imm\to\tZ_8$ is well defined by Lemma ~\ref{breg4}.

By the argument above, $\beta$ is additive under connected
summation. Hence, it is a homomorphism. 
\end{pf}

\subsection{A homomorphism $\Imm\to\tZ_4$}\label{tau}
Let $f\colon S^3\to\tR^5$ be a generic immersion.
Choose, once and for all, an isomorphism
$$
\psi\colon\pi_1(SO(5)/\tA(5))\to \tZ_4.
$$
(This choice is discussed in Section ~\ref{secrelinv}.)
\begin{dfn}
A component $c$ of the self intersection $M_f$
is {\em right twist framed} if 
$\psi(\omega_f(c))=1$. It is {\em left twist framed} if
$\psi(\omega_f(c))=-1$. 
\end{dfn}
If $c_S$ is disconnected then there is an induced spin structure on
$c$, see \cite{E1}, Proposition 2.10. This spin 
structure on $c$  is trivial (spin-cobordant to zero) if
$\omega_f(c)=0$ and nontrivial if 
$\omega_f(c)=p_\ast(u)$ (see \cite{E1}, proof of Theorem 7.30). 

We define $r(f)$, $l(f)$ and $n(f)$ as the 
number of self intersection components of a generic immersion $f$ that
are right twist framed, left twist framed and has the nontrivial spin
structure, respectively.
Define
$$
\tau(f)=r(f)+2n(f)-l(f)\mod{4}\in\tZ_4
$$
\begin{prp}\label{homtau}
For $\xi\in\Imm$, let $\tau(\xi)=\tau(f)$, where $f$ is a generic
immersion representing $\xi$. Then
$$
\tau\colon\Imm\to\tZ_4,
$$
is a homomorphism.
\end{prp}
\begin{pf}
That $\tau$ is well defined follows from \cite{E2}, Propositions
5.1.1-4. Clearly, $\tau$ is additive under connected sum. Hence, it 
is a homomorphism.
\end{pf}

\section{Special immersions}  
\noindent
In this section,
we construct immersions $S^3\to\tR^5$ with arbitrary Smale invariant
by rotating a 2-disk with a kink i 5-space and use an immersion
$\tR P^3\to\tR^4$ to show that the homomorphism $\beta$ (see Proposition
~\ref{hombeta}) is nontrivial. 

\subsection{The Whitney kink.}\label{whki}
\noindent
As in Whitney \cite{W2} we consider the map
$g\colon\tR^2\to\tR^4$ given by the equations
\begin{alignat*}{2}
g_1(x,y)&=x - \frac{2x}{u}  &\quad\quad  g_2(x,y)&=y\\
g_3(x,y)&=\frac{1}{u}  &  g_4(x,y)&=\frac{xy}{u},
\end{alignat*}
where $u=(1+x^2)(1+y^2)$.
This is an immersion with one transversal double point:
$$
(0,0,1,0)=g(1,0)=g(-1,0).
$$
Moreover, if $|x|$ or $|y|$ is large then $g$ is close to the standard
embedding $\tR^2\to\tR^2\times 0\subset\tR^4$. Thus, we can change
$g$ slightly so that it agrees with the standard embedding
outside some large disk. 

The Whitney kink enjoys the following symmetry property:
If $L\colon\tR^4\to\tR^4$ and $R\colon\tR^2\to\tR^2$ are the linear
maps 
$$
L(z_1,z_2,z_3,z_4)=(-z_1,-z_2,z_3,z_4)\quad\text{and}\quad R(x,y)=(-x,-y),
$$ 
then $L\circ g\circ R=g$. 

The differential of $g$ is 
$$
dg_{(x,y)}=
\begin{bmatrix}
1-\frac{2(1-x^2)}{(1+x^2)u} &
\frac{4xy}{(1+y^2)u}\\ \vspace{4pt}
0                           &        1            \\ \vspace{3pt}
\frac{-2x}{(1+x^2)u}        & \frac{-2y}{(1+y^2)u}\\ \vspace{3pt}
\frac{y(1-x^2)}{(1+x^2)u}   &
\frac{x(1-y^2)}{(1+y^2)u} 
\end{bmatrix}
$$
Let $\partial_x$ and $\partial_y$ be unit vector fields on $\tR^2$ in
the $x$ and $y$ directions, respectively. Applying $dg$ to these we get
a map of $\tR^2$ into $V_{4,2}$, which outside some large disk is constantly
equal to $(\partial_{1},\partial_{2})$, where $\partial_i$ is the unit
vector in the $z_i$-direction in $\tR^4$. 

Clearly, this map is homotopic in $V_{4,2}$ to a map of the form
$(v_1(x,y),\partial_{2})$, where $v_1(x,y)$ is orthonormal to
$\partial_{2}$. 
Let $D$ be a large disk in $\tR^2$. A straightforward calculation
shows that the degree of $v_1\colon (D^2,\partial
D^2)\to S^2$ is one.

Let $B$ be a 2-dimensional hemisphere in $\tR^4$ which is flat close to
its north pole. To make this precise, fix a small $a>0$ and let $B$ be a
2-disk such that 
\begin{align*}
&B\cap\left\{ z\in\tR^4 : z_3\le\sqrt{1-16a^2} \right\} =\\
&=\left\{ (\xi,\eta,\zeta,0) : \xi^2+\eta^2+\zeta^2=1, 0\le \zeta\le
\sqrt{1-16a^2} \right\},\\ 
&B\cap\left\{z\in\tR^4 : z_3\ge\sqrt{1-9a^2}\right\} =\\
&=\left\{(\xi,\eta,\sqrt{1-9a^2},0) : \xi^2+\eta^2\le 4a^2 \right\}.
\end{align*}
Let $(\xi,\eta)$, $\xi^2+\eta^2\le 1$, as above be coordinates on
$B$. Then we can define an 
immersion $K\colon B\to\tR^4$ with the following properties: The map
$K$ equals the inclusion for $\xi^2+\eta^2\ge a^2$. 
On the remaining part $D_a$ of $B$ the map $K$ is a suitably scaled
Whitney kink and $L(K(R(\xi,\eta)))=K(\xi,\eta)$. 

Let $(r,\vartheta)$, $0\le r\le\pi$ and $0\le\vartheta\le2\pi$ be
polar coordinates on the disk $D_a$, then the map $D_a\to V_{4,2}$
induced by the differential of $K$ is homotopic to
\begin{align*}
dK_{(r,\vartheta)}(\partial_{\xi}) &\simeq 
-\cos(r)\,\partial_{1} + \sin(r)\left(\cos
(\vartheta)\,\partial_{3}-\sin(\vartheta)\,\partial_{4}\right),\\
dK_{(r,\vartheta)}(\partial_{\eta}) &\simeq \partial_{2},
\tag{$\dagger$}
\end{align*}
since the first of these equations defines a map $(D_a,\partial
D_a)\to S^2$ of degree one.

\subsection{A modified standard embedding.}
\noindent
Let $(y_1,\dots,y_5)$ be coordinates on $\tR^5$. Let $\partial_i$
be the unit vector in the $y_i$-direction. Let $0\le\theta\le 2\pi$
and let $\tR^4_+(\theta)\subset\tR^5$ be the subset
$$
\{(x_0,x_1,x_2\cos\theta,x_2\sin\theta,x_3): x_2\ge 0\}.
$$ 
Let $D(\theta)\subset\tR_+^4(\theta)$ be the disk
$$
\{(x_0,x_1,x_2\cos\theta,x_2\sin\theta,0): x_2\ge 0, x_0^2+x_1^2+x_2^2=1\}.
$$
As $\theta$ varies from $0$ to $2\pi$, $D(\theta)$ sweeps the
standard $S^3$ in $\tR^5$. We note that if $\theta_1\ne\theta_2$,
$\theta_1\ne0$  then
$$
D(\theta_1)\cap D(\theta_2)=\{(x_0,x_1,0,0,0): x_0^2+x_1^2=1\},
$$
and $D(0)=D(2\pi)$. 

A straightforward calculation shows that the column vectors
$S^i(\theta)$ of
the matrix
$S(0,0,\cos\theta,\sin\theta,0,0)=\Theta^s(0,0,\cos\theta,\sin\theta,0,0)$
(see Lemma ~\ref{calcSmin5} for notation)
are  
\begin{align*}
S^1(\theta) &= \sin(\theta)\,\partial_3 -\cos(\theta)\,\partial_4,\\
S^2(\theta) &= -\cos(\theta)\,\partial_1 -\sin(\theta)\,\partial_2,\\
S^3(\theta) &= -\sin(\theta)\,\partial_1 +\cos(\theta)\,\partial_3,\\
S^4(\theta) &= \partial_5,\\
S^5(\theta) &= \cos(\theta)\partial_3 +\sin(\theta)\partial_4.
\end{align*}

We now modify the standard embedding: Replace the disk
$D(\theta)$ by $B(\theta)$, where $B(\theta)$ is obtained by
rotating $B\subset\tR^4_+(0)$ to $\tR^4_+(\theta)$. The embedding of
$S^3$ so obtained is clearly diffeotopic to the standard embedding. 
We may assume that the framing it induces on $D_a(\theta)$ depends only
on $\theta$ and equals $S(\theta)$.

\subsection{Constructing immersions with arbitrary Smale
invariant}\label{secgen} 
\noindent
Let the linear map 
$L_\alpha\colon\tR^5\to\tR^5$ be given by the matrix
$$
\begin{bmatrix}
\begin{smallmatrix}\cos\alpha & -\sin\alpha\\ 
\sin\alpha & \cos\alpha \end{smallmatrix} & 
\begin{smallmatrix}{} & {} &{}\\ 
{} & {} &{} \end{smallmatrix}\\  
\begin{smallmatrix}{} & {} \\ 
{} & {} \\ 
{} & {} \end{smallmatrix}  
&\tI_3
\end{bmatrix}.
$$
Let $(\xi,\eta)$ be coordinates on $B(\theta)$, as in
Section ~\ref{whki}.  
Let $R_\alpha\colon B(\theta)\to B(\theta)$ be the map that is
given by the matrix 
$$
\begin{bmatrix}\cos\alpha & -\sin\alpha\\ 
\sin\alpha & \cos\alpha \end{bmatrix},
$$
in the coordinates $(\xi,\eta)$.

Consider $S^3=\bigcup_{0\le\theta\le2\pi} B(\theta)$. 
For every half integer $m\in\tfrac12\tZ$ we define $f_m\colon S^3\to\tR^5$,
$$
f_m|B(\theta)=L_{m\theta}\circ K\circ R_{-m\theta},
$$
where $K$ is as in Section \ref{whki}.
Then $f_m$ is well defined since it agrees with the inclusion on
$\partial B(\theta)$ for each $\theta$ and since $K$ is symmetric, i.e.
$L_{2m\pi}\circ K\circ R_{-2m\pi}=K$.

The self intersection of $f_m$ is one circle. It is traced
out by the double point of $K$ under rotation. The
preimage of this circle is a torus knot of type $(2m,2)$ if $m$ is not an
integer and a two component link with unknotted components that are
linked with linking number $m$ if $m$ is an integer.
\begin{prp}\label{theg}
The Smale invariant $\Omega(f_m)$ satisfies,
$$
\Omega(f_m)=-2m[\sigma]+N\in\pi_3(V_{5,3}).
$$
(See Lemma ~\ref{pi3V53} for notation.)
\end{prp}
\noindent
Thus, as $m$ runs through $\tfrac12\tZ$, $f_m$ runs through the regular
homotopy classes of immersions $S^3\to\tR^5$.
The proof of Proposition ~\ref{theg}
constitutes Section ~\ref{pftheg}. 
\begin{cor}
The homomorphism (see Proposition ~\ref{homtau})
$$
\tau\colon\Imm\to\tZ_4,
$$
is surjective.
\end{cor}
\begin{pf}
The immersion $f_{\frac12}$ represents a generator of $\Imm$. It has
one self intersection circle with connected preimage. Thus,
$\tau(f_{\frac12})$ is a generator of $\tZ_4$.
\end{pf}
\begin{lma}\label{lkgen} 
The invariant $\lk$ (see Lemma ~\ref{lk}) takes the values
$$
\lk(f_m)=-4m.
$$
\end{lma}
We prove Lemma ~\ref{lkgen} in Section ~\ref{pflkgen}.
\begin{cor}
The homomorphism (see Corollary ~\ref{homlambda})
$$
\lambda\colon\Imm\to\tZ_3,
$$
is surjective.
\end{cor}
\begin{pf}
By Lemma ~\ref{lkgen}, $\lambda(f_{\frac12})$ is a generator of $\tZ_3$.
\end{pf}

\subsection{An immersion into 4-space}
\noindent
As in Lashof and Smale \cite{LS}, Theorem 3.4, we 
consider an immersion $k\colon S^2\to\tR^4$ with one transversal double
point. The normal bundle of this immersion has Euler number 2 and
the boundary of a small
tubular neighborhood of it is therefore an immersion
$h\colon US^2\cong\tR P^3\to \tR^4$, where $US^2$ denotes the unit
tangent bundle of $S^2$. 

Let $\pi\colon S^3\to\tR P^3$ be the universal double covering. 
Let $f=h\circ\pi\colon S^3\to\tR^4$. Then $f$ is an immersion. 
Let $g\colon S^3\to\tR^4$ be a generic immersion regularly homotopic
to $f$. 
\begin{prp}\label{RP3} 
The invariant $\beta(g)$ (see Proposition ~\ref{breg4}) is a
generator of $\tZ_8$, the immersion $g$  
has an odd number of 
quadruple points and the Euler characteristic of $F_g$ is odd.
\end{prp}
\noindent
Proposition ~\ref{RP3} is proved in Section ~\ref{pfRP3}
\begin{cor}
The $\tZ_2$-valued regular homotopy invariants $Q_2$ and $D_2$ defined
in Section ~\ref{secmorinv} are non-trivial.\qed
\end{cor}
\begin{cor}
The homomorphism (see Proposition ~\ref{hombeta})
$$
\beta\colon\Imm\to\tZ_8
$$
is an epimorphism.\qed
\end{cor}
\section{Proofs}
\noindent
In this section we prove Theorems
~\ref{thmimem}, ~\ref{gi41},~\ref{gi42}, and ~\ref{gi43}, Proposition 
~\ref{theg}, Lemma ~\ref{lkgen}, and Proposition ~\ref{RP3}

\subsection{Proof of Theorem ~\ref{thmimem}}\label{pfimem}
\noindent
Proposition ~\ref{sgntr} shows that $\sigma$ is an
isomorphism. Proposition ~\ref{24} shows that $\Emb$ is a subgroup of
the infinite cyclic group $\Imm$ of index 24. The Smale invariant $\Omega$
takes values in $\pi_3(V_{5,3})$. Choosing a generator, we identify
this group with $\tZ$. For one of the two possible choices, the first
part of the diagram is  
correct. (Making the other choice, the second homomorphism 
in the second row would be $\times(-24)$ instead of $\times 24$). 

To prove the theorem it is enough to show that $\tau\oplus\beta$ is
surjective. 
Lemma ~\ref{lkgen} shows that
$\tau(f_{\frac12})$ is a generator of $\tZ_3$ and if $g\colon
S^3\to\tR^4\subset\tR^5$ is as in Proposition
~\ref{RP3} then $\beta(g)$ is a generator of $\tZ_8$. 
Hence, if either $\beta(f_{\frac12})$ is a generator of 
$\tZ_8$ or $\tau(g)$ is a generator of $\tZ_3$ then $\tau\oplus\beta$ is
surjective. Assume that this is not the case. Then
$\beta(f_{\frac12})$ is even and 
$\tau(g)=0$ and hence 
$$
(\tau(f_{\frac12}\star g),\beta(f_{\frac12}\star g))=
(\tau(f_{\frac12}),\beta(f_{\frac12}))+
(\tau(g),\beta(g)),
$$ 
is a
generator of $\tZ_3\oplus\tZ_8$. Thus, $\tau\oplus\beta$ is
surjective. (Since $f_{\frac12}$ represents a generator of $\Imm$, it
actually follows that $\beta(f_{\frac12})$ is odd.)
\qed
\subsection{Proofs of Theorems ~\ref{gi41},~\ref{gi42}, and
~\ref{gi43}}\label{secrelinv} 
\noindent
Consider the invariants $\beta$, $\tau$,
$D_2$ and $Q_2$, defined in Section ~\ref{secmorinv}.

Since $D_2$ and $Q_2$ are both additive under connected summation (see
the proof of Lemma ~\ref{hombeta}), it
follows exactly as for $\beta$ that they are invariant under regular
homotopy in $\tR^5$. They both define homomorphisms $\Imm\to\tZ_2$ and
Lemma ~\ref{RP3} implies that they are both nontrivial and therefore
equal.  

{\em Proof of Theorem ~\ref{gi41}:} Let $f\colon S^3\to\tR^4$ be a
generic immersion. Consider $f$ as an immersion into $\tR^5$ and
evaluate the $\tZ_2$-valued invariants induced by $D_2$ and $Q_2$. By
the above they are equal. The theorem follows.\qed

Let $\mu\colon\Imm\to\tZ_2$ denote the homomorphism obtained from $Q_2$
and $D_2$. We have the following sequence:
$$
\Imm\xrightarrow{\beta}\tZ_8\xrightarrow{r_{4}}\tZ_4\xrightarrow{r_{2}}\tZ_2,
$$
where $r_2$ and $r_4$ denote reduction modulo 2 and reduction modulo
4, respectively. 
The Brown invariant of a pin structure on a surface reduced modulo 2
equals the Euler characteristic of the surface reduced modulo 2.
Therefore, $r_2\circ r_4\circ\beta=\mu$. 

Moreover, if $f\colon S^3\to\tR^5$ represents a generator of $\Imm$
then $r_4(\beta(f))$ is a generator of $\tZ_4$ and so is $\tau(f)$. 

As a consequence,  
$$
\beta(h)\mod{4}=\pm \tau(h)\in\tZ_4
$$ 
for any immersion $h\colon S^3\to\tR^5$. The sign in this formula
depends on the choices of the isomorphisms $\phi$ (Section
~\ref{secbeta}) and $\psi$ (Section ~\ref{tau}). For any one of the
two possible choices of $\phi$ there is a unique choice (out of two
possible) of $\psi$ such that the sign in the above formula is $+$.

{\em Proof of Theorem ~\ref{gi42}:} Let $f\colon S^3\to\tR^5$ be a
generic immersion and $g\colon S^3\to\tR^4$ a generic immersion
regularly homotopic to $f$. Then $\tau(f)$ is odd if and only if
$r_4(\beta(g))$ is odd. But $r_4(\beta(g))$ is odd if and only if
$\beta(g)$ is odd which is equivalent to $F_g$ having odd Euler
characteristic. Apply Theorem ~\ref{gi41}.\qed

{\em Proof of Theorem ~\ref{gi43}:} Let $f$ and $g$ be as above.
If $F_g$ is orientable then $\beta(g)$ is divisible by $4$. Hence,
$r_4(\beta)=0$ and therefore $\tau(f)=0$.\qed  

\subsection{Proof of Proposition \ref{theg}}\label{pftheg}
\noindent
For simplicity of notation we let $s\colon S^3\to \tR^5$ denote the
{\em modified}
standard embedding and $S\colon S^3\to SO(5)$ denote the map $\Theta^s$
associated to it. Fix $m$ and let
$F\colon S^3\to SO(5)$ denote the 
map $\Theta^{f_m}$ associated to the immersion $f_m\colon
S^3\to\tR^5$. Let $\Lambda_m(x)=S^{-1}(x)F(x)$, where $S^{-1}(x)$ is
the inverse of the matrix $S(x)$ and $S^{-1}(x)F(x)$ is the matrix product.
To compute the Smale invariant of $f_m$ we must compute the homotopy
class (see Lemma ~\ref{calcSmin5})
$$
[\Theta^{f_m}]-[\Theta^s]=[F]-[S]=[\Lambda_m]\in\pi_3(SO(5))=\pi_3(V_{5,3}).
$$
Since $f_m=s$ outside the solid torus
$T=\bigcup_{0\le\theta\le2\pi} D_a(\theta)$, we have
$$
\Lambda_m(x)=S^{-1}(x)F(x)=
\begin{bmatrix}
\tI_3 & 
\begin{smallmatrix}
{}&{}\\
{}&{}\\
{}&{}
\end{smallmatrix}\\
\begin{smallmatrix}
{}&{}&{}\\
{}&{}&{}
\end{smallmatrix}
 & \tJ(x)
\end{bmatrix},\text{ for }x\in S^3-T,
$$
where $\tJ(x)\in SO(2)$. 

Let $(\theta,r,\varphi)$,
$\theta,\varphi\in[0,2\pi]$, $r\in[0,\pi]$ be coordinates on
$T$. Here $\theta$ indicates which $D_a(\theta)$ we are in and on
$D_a(\theta)$ we have
$$
\xi=\tfrac{ra}{\pi}\cos(\varphi),\quad
\eta=\tfrac{ra}{\pi}\sin(\varphi).
$$
Using equations ($\dagger$) in Section ~\ref{whki} we can compute the
first columns of 
$F|T$. The result is the following:
\begin{align*}
F^1(\theta,r,\varphi) =&
\sin(\theta)\,\partial_3 -\cos(\theta)\,\partial_4,\\
F^2(\theta,r,\varphi) =&-\cos(\theta-m\theta))\,v_1(\theta,r,\varphi)
-\sin (\theta-m\theta)\,v_2(\theta,r,\varphi)\\
F^3(\theta,r,\varphi) =&-\sin(\theta-m\theta)\,v_1(\theta,r,\varphi)
+\cos(\theta-m\theta)\,v_2(\theta,r,\varphi),
\end{align*}
where
\begin{align*}
v_1(\theta,r,\phi) =& -\cos(r)
\bigl[\cos (m\theta)\,\partial_1+\sin (m\theta)\,\partial_2\bigr]\\
&+\sin(r)
\bigl[\cos(\varphi-m\theta)
\left(\cos (\theta)\,\partial_3+\sin (\theta)\,\partial_4\right)
-\sin(\varphi-m\theta)\,\partial_5
\bigr],\\
v_2(\theta,r,\varphi)=& -\sin(m\theta)\,\partial_1+\cos(m\theta)\,\partial_2.
\end{align*}
We note that $F^1(x)=S^1(x)$, for all $x\in S^3$. Hence, 
$$
\Lambda_m(x)=
\begin{bmatrix}
1 & \quad \\
\quad & \tM(x)
\end{bmatrix},\quad\text{for all }x\in S^3,
$$
where $\tM(x)\in SO(4)$. Thus, to compute $[\Lambda_m]\in\pi_3(SO(5))$ it
is enough to compute the homotopy class of $\tM\colon S^3\to SO(4)$ which,  
as we shall see, can be calculated if we know the first two columns of
$\tM(x)$. 

If $\tM=(m^j_i)$ then $m_j^i=\langle S^{i-1},F^{j-1}\rangle$,
$i,j=2,3,4,5$, where $\langle{ },{ }\rangle$ is the usual inner product
on $\tR^5$. Hence, we have all the information we need to determine
$[\Lambda_m]$. In what follows, we first study how to calculate the
homotopy 
class of $\tM$ and then carry out the necessary calculations.
 
In terms of our standard generators of $\pi_3(SO(4))$ we have, for
$\tM\colon S^3\to SO(4)$, $[\tM]=u[\sigma]+v[\rho]\in\pi_3(SO(4))$,
$u,v\in\tZ$. 

The integer $u$ is simply the degree of the map $p\circ\tM=\tM^1\colon
S^3\to S^3$, where $SO(4)\xrightarrow{p} S^3$ is the fibration
described in Section ~\ref{sechomth} and $\tM^1(x)$ is the first column vector of
$\tM(x)$. 

To compute $v[\rho]$, recall from Section ~\ref{sechomth} that
$SO(4)\cong S^3\times SO(3)$ and thus, we must compare the homotopy
class of the map 
$p_2\circ\tM\colon S^3\to SO(3)$, where $p_2\colon SO(4)\to SO(3)$ is
the projection onto the second factor in the product space $SO(4)\cong
S^3\times SO(3)$, to that of $p_2\rho\colon S^3\to SO(3)$. 

The product structure on $SO(4)$ is obtained by using the section
$\sigma\colon S^3\to SO(4)$. The 
projection $p_2\colon SO(4)\to SO(3)$ is thus determined by the
equation  
$$
\sigma(p(x))^{-1}x=
\begin{bmatrix}
1 & \quad \\
\quad & p_2(x)
\end{bmatrix},\quad x\in SO(4).
$$
We note that $p_2\circ\rho=\varrho\colon S^3\to SO(3)$. Let
$\tN=p_2\circ\tM$. 
To compare the homotopy classes of these maps we can use the fibration
$SO(3)\stackrel{p}{\to}S^2$: Since
$p_\ast\colon\pi_3(SO(3))\to\pi_3(S^2)$ is an isomorphism we might as
well compare the homotopy classes of the maps $p\circ\varrho$ and
$p\circ\tN$ both mapping $S^3$ to $S^2$. By the Pontryagin
construction the homotopy class of any map $S^3\to S^2$ is determined
by its Hopf invariant (the linking number of two regular fibers).

The map $\varrho\colon S^3\to SO(3)$,
has Hopf invariant one if endow $S^3$ with the orientation that is
coherent with the quaternion framing.

Thus, to determine $[\tN]\in\pi_3(SO(3))$ we need only calculate the
Hopf invariant of $p_2\circ\tN=\tN^1\colon S^3\to S^2$, where
$\tN^1(x)$ is the first column vector in the matrix $\tN(x)$. 
This, in turn, means that the integer $v$ can be computed as the Hopf
invariant of the map $\tN^1\colon S^3\to S^2$, where $S^3$ has the
same orientation as above.

We now turn to the actual calculations. 
We have $\tM^k=\sum_{i=1}^4 \langle
S^{i+1},F^{k+1}\rangle\,\partial_i$, $k=1,2$ where $\partial_i$,
$i=1,\dots,4$ is the standard basis in $\tR^4$ and $\langle\,
,\,\rangle$ is the standard inner product on $\tR^5$. Computing the
necessary inner products we find that if $x\in T\subset S^3$ then
$$
\tM^1(x)=\tM^1(\theta,r,\varphi)=
\begin{bmatrix}
-\cos^2(m\theta-\theta)\cos(r)+\sin^2(m\theta-\theta)\\ \vspace{6pt}
\frac12\left(\cos(r)+1\right)\sin(2m\theta-2\theta)\\ \vspace{6pt}
\cos(m\theta-\theta)\sin(r)\sin(\varphi-m\theta)\\ \vspace{6pt}
-\cos(m\theta-\theta)\sin(r)\cos(\varphi-m\theta)
\end{bmatrix},
$$
$$
\tM^2(x)=\tM^2(\theta,r,\varphi)=
\begin{bmatrix}
\frac12\left(\cos(r)+1\right)\sin(2m\theta-2\theta)\\ \vspace{6pt}
-\sin^2(m\theta-\theta)\cos(r)+\cos^2(m\theta-\theta)\\ \vspace{6pt}
-\sin(m\theta-\theta)\sin(r)\sin(\varphi-m\theta)\\ \vspace{6pt}
\sin(m\theta-\theta)\sin(r)\cos(\varphi-m\theta)
\end{bmatrix}
$$
and if $x\in S^3-T$ then
$$
\tM^1(x)=\begin{bmatrix}
1\\ 0\\ 0\\ 0\\
\end{bmatrix},\quad 
\tM^2(x)=\begin{bmatrix}
0\\ 1\\ 0\\ 0\\
\end{bmatrix}.
$$
To find the first column of $\tN$ we must calculate
$\sigma(p(\tM))^{-1}\tM^2$. For $y\in\tR^4$, $\sigma(p(\tM))y=\tM^1\cdot y$,
(for notation see Section ~\ref{sechomth}). For $q\in\tH$, let $\overline
q$ denote its conjugate then  $\sigma(p(\tM))^{-1}y=\overline{\tM^1}\cdot
y$. Multiplying gives
$$
\overline{\tM^1(x)}\cdot\tM^2(x)=
\begin{bmatrix}
0\\
-\cos(r)\\ \vspace{3pt}
\sin(r)\cos(\varphi-\theta)\\ \vspace{3pt}
\sin(r)\sin(\varphi-\theta)
\end{bmatrix}\text{ if $x\in T$.}
$$
Hence, the the first column of $\tN$ is given by
$$
\tN^1(x)=
\begin{bmatrix}
-\cos(r)\\ \vspace{3pt}
\sin(r)\cos(\varphi-\theta)\\ \vspace{3pt}
\sin(r)\sin(\varphi-\theta)
\end{bmatrix}\text{ if $x\in T$ and }
\tN^1(x)=
\begin{bmatrix}
1\\
0\\
0\\
\end{bmatrix}\text{ if $x\in S^3-T$.}
$$

To determine the degree of $\tM^1$ we consider preimages of the regular
value $(0,0,1,0)\in S^3$. For $x\in S^3-T$, $\tM^1(x)\ne (0,0,1,0)$
and for $x\in T$ we get the equations 
\begin{align*}
r&=\frac{\pi}{2},\\
\theta&=\frac{l\pi}{m-1},\\ 
\varphi&=\frac{ml\pi}{m-1}+(-1)^{l}\,\frac{\pi}{2},
\end{align*}
where $l= 0,1,\dots,|2m-3|$. Thus, there are $|2m-2|$ points in the
preimage of $(0,0,1,0)$.

The orientation of $T$ induced from the quaternion framing on $S^3$ is
the same as that given by the basis
$(\partial_\theta,\partial_r,\partial_\varphi)$. Thus, the orientation of
$S^3$ at $(0,0,1,0)$ induced by $\tM^1$ is
$o((m-1)\,\partial_2,\partial_1,\partial_4)$. The orientation of $S^3$
at $(0,0,1,0)$ induced from the standard orientation on $\tR^4$ is
$o(\partial_1,\partial_2,\partial_4)$. So the sign of each point in
the preimage is the same. It is positive if $m<1$ and negative if
$m>1$. Thus, $\deg(\tM^1)=u=-(2m-2)=-2m+2$.  

To compute the Hopf invariant of $\tN^1\colon S^3\to S^2$ we
consider the preimage of the regular value $(0,0,1)$ and see that it
is the curve $c$ defined by the equations 
$$
r=\frac{\pi}{2},\quad\quad \phi-\theta=\frac{\pi}{2}.
$$
Moreover, the vector field $-\frac{\partial}{\partial r}$ along $c$, is
mapped by $d\tN$ to the vector $(1,0,0)$. Shifting the curve $c$ along this
vector field we get a curve $c'$ and $\lk(c,c')=-1$ using the above
orientation on $T$. So, the Hopf invariant of $\tN^1$ is $-1$ and thus
$v=-1$.  

Collecting these results, we finally get
\begin{align*}
\Omega(f_m) &= [\Lambda_m] = [\tM]+N = ((-2m+2)[\sigma]-[\rho])+N\\
&=-2m[\sigma]+N\in\pi_3(SO(5))=\pi_3(V_{5,3}),
\end{align*}
as claimed.\qed
\subsection{Proof of Lemma ~\ref{lkgen}}\label{pflkgen}
Let $(r,\theta,\varphi)$ be coordinates on the solid torus where the
Whitney kinks are placed. The preimage of the self intersection $c_m$
of $f_m$ is
$r=\epsilon$, $\varphi=m\theta$ or $\varphi=m\theta+\pi$, where
$\epsilon$ is some small positive number. Consider the vector field
$\partial_r$ along this preimage. If we shift $c_m$ along 
$$ 
v(p)=(df_m(\partial_r(p_1))+df_m(\partial_r(p_2)),
$$
where $f_m^{-1}(p)=\{p_1,p_2\}$ then we obtain a curve $c''_m$, which
bounds a 2-disk in the $x_3x_4$-plane without intersections with
$f_m(S^3)$. Thus, $\lk(c''_m,f_m(S^3))=0$. 

Assume that $m$ is not an integer. 
Then $c_S=f_m^{-1}(c_m)$ is connected. 
Shifting $c_S$ along $\partial_r$ gives a curve $c_S''$ with
$\lk(c_S,c_S'')=4m$. Thus, the framing of $c_S$ which satisfies
condition $(\dagger)$ in Section ~\ref{seclk} is the one that differs
from $\partial_r$ by $-4m$ rotations. 

Assume that $m$ is an integer. Then $c_S=c_1\cup c_2$ have two
components and $\lk(c_1,c_2)=m$. Shifting $c_i$ along $\partial_r$
gives a curve $c_i''$ with $\lk(c_i,c_i'')=m$. Thus, the framings of
the $c_i$ that satisfy condition $(\dagger)$ in Section ~\ref{seclk} differ
from $\partial_r$ by $-2m$ rotations each.

The rotations that we need to add to $\partial_r$ can be made locally
supported. Using a local model it is easy to see how they affect the
linking number:
Let $X$ and $Y$ be the 3-planes $X=(x_1,x_2,x_3,0,0)$ and
$Y=(y_1,0,0,y_2,y_3)$. Then $X\cap Y$ is the line
$l=(t,0,0,0,0)$. Shifting $l$ along
$\epsilon(\partial_2+\partial_4)$ we get
$l'=(t,\epsilon,0,\epsilon,0)$. Shifting $l$ along
$\epsilon(-\cos(t)\partial_2+\sin(t)\partial_3+\partial_4)$, $-\pi\le
t\le\pi$ we get the 
curve $l''=(t,-\epsilon\cos(t),\epsilon\sin(t),0,\epsilon,0)$. Let 
$D=(1-s)l'+sl''$, $0\le s\le 1$. Then $D$ is a 2-chain with boundary
$l''-l'$, $D\cap X=\emptyset$ and $D\cap Y=(0,0,0,\epsilon,0)$ with
sign $-1$. Thus, for each rotation in the negative direction we
decrease the linking number by 1.

Adding the appropriate number of rotations to the framing $\partial_r$
and shifting the curve $c_m$ along the corresponding vector field we
obtain a curve $c'_m$ with
$$
\lk(c_m',f_m(S^3))=\lk(c''_m,f_m(S^3))-4m=-4m.
$$
Thus, $\lk(f_m)=-4m$.\qed

\subsection{Proof of Proposition ~\ref{RP3}}\label{pfRP3}
\noindent
The self intersection of $h$ is located close to the transverse double
point $h(p)=h(q)=0$ of $S^2$. Let $x$ and $y$ be coordinates on small
disks $D_1$ and $D_2$ around
$p$ and $q$ respectively. We may assume that 
$$
k(x)=(x_1,x_2,0,0)\quad\text{and}\quad k(y)=(0,0,y_1,y_2).
$$
Close to $p$ and $q$, $US^2$ is of the form $D_i^2\times S^1$, $i=1,2$ and the
immersion $h$  is:
$$
h(x,\theta)=(x_1,x_2,\epsilon\cos(\theta),\epsilon\sin(\theta)),\quad\quad
h(y,\varphi)=(\epsilon\cos(\varphi),\epsilon\sin(\varphi),y_1,y_2),
$$
where $0\le\theta\le2\pi$ and $0\le\varphi\le2\pi$ are coordinates on
$S^1$ in the fiber close to $p$ and $q$, respectively.

Thus, $M_h=T^2$ and the preimage ${\wt M}_h$ is $c_i\times S^1\subset
D_i\times S^1$, where $c_i$ is a circle of radius $\epsilon$ around
$0\in D_i$.  

We can find $\tR P^2\subset US^2$ as follows: Fix a 
point $s$ in $S^2$ and a unit tangent vector $v$ at that point. 
The space of 2-planes $\Pi$ through $0$ in $\tR^3$ is $\tR P^2$. 
Let $r_\Pi$ be the reflection in $\Pi$. The map
$\Pi\mapsto(r_\Pi(s),r_\Pi(v))$ is an embedding of $\tR P^2$ into $US^2$.

Consider the double cover $p\colon S^3\to US^2$. Observe that
$K=p^{-1}(\tR P^2)\cong S^2$ subdivides $S^3$ into
two hemispheres $D^3_+$ and $D^3_-$.

The intersection ${\wt M}_h\cap \tR P^2$ consists of two circles,
isotopic to $c_i\times{1}\subset D_i\times S^1$, $i=1,2$. Altering the
embedding of $\tR P^2$ slightly we may assume that ${\wt M}_h\cap\tR
P^2=c_1\times 1\cup c_2\times 1$.  

The covering $p\colon S^3\to US^2$ is nontrivial on fibers so 
$p^{-1}({\wt M_h})$ consists of two tori $T_i\subset
S^3$, $i=1,2$ and $p|T_i\colon T_i\to c_i\times S^1$ is a double cover. Thus,
$K\cap T_i$  consists of two meridians in $T_i$, $i=1,2$.  

Choose a normal vector field $\nu$ to $h$. Let $\delta>0$ be small and
let $\phi\colon S^3\to(-\delta,\delta)$ be positive on $D_+$, negative
on $D_-$, $0$ on $K$, and have nonvanishing derivative in
the normal direction of $K$. (Map $(S^3,K)$ to
$(S^3,\{x_4=0\})\subset\tR^4$ and use a suitably scaled height
function.) We use $\phi$ and $\nu$ to perturb $f$. For $x\in S^3$, let
$$
g(x)=f(x)+\phi(x)\nu(x).
$$
Then $g$ is an immersion. Clearly, $g$ is an embedding outside a small
neighborhood of $T_1\cup T_2\cup K$ in $S^3$. Moreover, inside this
neighborhood it has only transverse double points if we stay away from
$T_1\cup T_2$. To determine the self intersection of $g$ we can
therefore restrict attention to neighborhoods of $T_1$ and $T_2$. Let
$E_i\times S^1=p^{-1}(D_i\times S^1)$. If $(x_1,x_2,e^{i\alpha})$,
$(y_1,y_2,e^{i\omega})$ are coordinates on $E_1\times S^1$ and
$E_2\times S^1$, respectively then $p\colon E_i\times S^1\to D_i\times
S^1$ is given by 
\begin{align*}
p(x_1,x_2,e^{i\alpha})&=(x_1,x_2,e^{i2\alpha}),\\
p(y_1,y_2,e^{i\omega})&=(y_1,y_2,e^{i2\omega}).
\end{align*}
Normal fields along $h(D_1\times S^1)$ and $h(D_2\times S^1)$ are
\begin{align*}
\nu_1(x_1,x_2,e^{i\theta})
&=\cos(\theta)\,\partial_3+\sin(\theta)\,\partial_4,\\ 
\nu_1(y_1,y_2,e^{i\varphi})
&=\cos(\theta)\,\partial_1+\sin(\theta)\,\partial_2. 
\end{align*}
We may assume that the function $\phi|E_i\times S^1$ is given by 
$$
\phi(x_1,x_2,e^{i\alpha})=\delta\sin(\alpha)\quad\text{and}\quad
\phi(y_1,y_2,e^{i\omega})=\delta\sin(\omega).
$$
Then, in $E_i\times S^1$, $g=f+\phi\nu$ is given by
\begin{align*}
g(x_1,x_2,e^{i\alpha})&=
\left(x_1,x_2,\left(\epsilon +\delta\sin(\alpha)\right)\cos(2\alpha),
\left(\epsilon +\delta\sin(\alpha)\right)\sin(2\alpha)\right),\\
g(y_1,y_2,e^{i\omega})&=
\left(\left(\epsilon +\delta\sin(\omega)\right)\cos(2\omega),
\left(\epsilon +\delta\sin(\omega)\right)\sin(2\omega),y_1,y_2\right).
\end{align*}
Thus, $g$ has one quadruple point at $(\epsilon,0,\epsilon,0)$, with
preimages 
$$
(x_1,x_2,e^{i\alpha})=(\epsilon,0,\pm1)=(y_1,y_2,e^{i\omega}),
$$
and $g$ has four triple lines:
\begin{align*}
\alpha &\mapsto
\left(\epsilon,0,
\left(\epsilon +\delta\sin(\alpha)\right)\cos(2\alpha),
\left(\epsilon +\delta\sin(\alpha)\right)\sin(2\alpha)\right), \\
\omega &\mapsto
\left(\left(\epsilon +\delta\sin(\omega)\right)\cos(2\omega),
\left(\epsilon +\delta\sin(\omega)\right)\sin(2\omega),
\epsilon,0\right),
\end{align*} 
where $\alpha,\omega\in(0,\pi)$ or $\alpha,\omega\in (\pi, 2\pi)$.

\begin{figure}[htbp]
\begin{center}
\includegraphics[angle=-90, width=10cm]{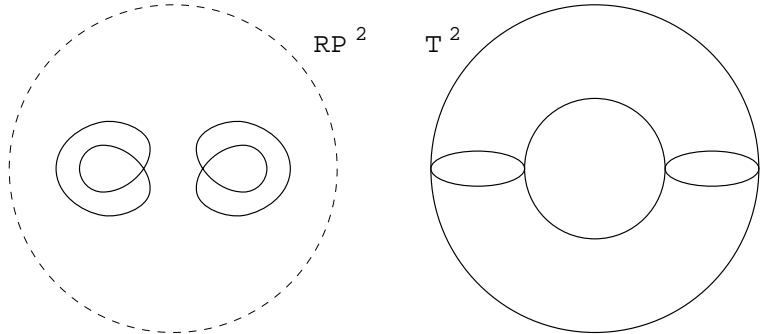}
\end{center}
\caption{The self intersection surface $F_g$}\label{figFg}
\end{figure}

Constructing $F_g$ by gluing the nine double point sheets of $g$
together gives $F_g=T^2\cup\tR P^2$ as shown in Figure ~\ref{figFg},
where the six 
intersection points of the solid lines represents $F_g^0$ and the
remaining twelve open arcs of the solid lines represents $F_g^1$.
It follows that $\beta(g)$ as well as $\chi(F_g)$ are odd. Since
$\beta$ and $D_2$ are invariants of regular homotopy in $\tR^5$, the
proposition follows. \qed

\section*{Acknowledgments}
\noindent
I want to thank Oleg Viro for valuable ideas and Ryszard Rubinsztein for
fruitful discussions. I also want to thank the referee for pointing
out that the problem on representing regular homotopy classes by
embeddings was solved in \cite{HM}. This shortened the paper considerably.


\begin{thebibliography}{99}
\bibitem{A1} V. I. Arnold, {\em The Vassiliev theory of discriminants
of knots},   
in First European Congress of Mathematics: Paris July 6-10,
1992, Volume I, Birkh{\"a}user (1994), 3-30.
\bibitem{A2} V. I. Arnold, {\em Plane curves, their invariants,
perestroikas and classifications},  
in Singularities and Bifurcations, edited by Arnold, Adv. Sov. Math. {\bf 21}
(1994) 39-91.
\bibitem{E1} T. Ekholm, {\em Immersions in the Metastable Range and Spin
Structures on Surfaces}, Math. Scand. {\bf 83} (1998) 5-41. 
\bibitem{E2} T. Ekholm, {\em Vassiliev invariants and regular homotopy
of generic immersions $S^k\to\tR^{2k-1}$, $k\ge 4$}, 
J. Knot Theor. Ramif. {\bf 7} (1998) 1041-1064.
\bibitem{E3} T. Ekholm, {\em Geometry of self intersection surfaces
and Vassiliev invariants of generic immersions $S^k\to\tR^{2k-2}$}, 
Preprint, Uppsala  University (1998)
\bibitem{Hi} M. W. Hirsch, {\em Immersions of manifolds},
Trans. Amer. math. Soc. {\bf 93} (1959), 242-276.
\bibitem{HM} J. F. Hughes and P. M. Melvin, {\em The Smale invariant of a
knot}, Comment. Math. Helv. {\bf 60} (1985), 615-627.
\bibitem{Hu} D. Husemoller, {\em Fibre Bundles, 3rd ed.}, Springer
Verlag (1994) 
\bibitem{K} M. A. Kervaire, {\em Sur le fibr\'e normal \`a une
sph\`ere immerg\'ee dans une espace euclidien},
Comment. Math. Helv. {\bf 33}, (1959) 121-131
\bibitem{KT} R. C. Kirby and L. R. Taylor, {\em {\em Pin} structures on 
low-dimensional manifolds}, in "Geometry of low-dimensional manifolds: 
2", London Math.\ Soc.\ Lect.\ Notes Ser. {\bf 151}, Cambridge University 
Press (1990), 177-242.
\bibitem{LS} R. Lashof and S.Smale, {\em On the immersion of manifolds
in Euclidean space}, Ann. of Math. {\bf 68} (1958) 562-583  
\bibitem{M} Y. Matsumoto, {\em An elementary proof of Rohlin's signature 
theorem and its extension by Guillou and Marin}, in "A la Recherche de la 
Topologie Perdue", edited by Guillou and Marin, Birkh{\"a}user (1986),
119-139.   
\bibitem{MK} J. W. Milnor and M. A. Kervaire, {\em Bernoulli numbers,
homotopy groups and a Theorem of Rohlin},
Proc. Int. Math. Congress. Edinburgh (1958), 454-458 
\bibitem{MS} J. W. Milnor and J. A. Stasheff, {\em Characteristic
Classes}, Ann.\ of Math.\ study {\bf 76}, Princeton University Press (1974)
\bibitem{S} S. Smale, {\em Classification of immersions of spheres in
Euclidean space},
Ann.\ of Math.\ {\bf 69} (1959), 327-344.
\bibitem{W2} H. Whitney, {\em The self intersections of a smooth
$n$-manifold in $2n$-space}, 
Ann. of Math. {\bf 45} (1944), 220-246.
\end{thebibliography}
\end{document}